
\input amstex
\documentstyle{amsppt}

\baselineskip=18truept
\voffset=0.3truein      
\hoffset=4truemm
\hsize=6.0truein        
\vsize=8.5truein        
\parskip=5truept        


\TagsOnRight

\def\today{\ifcase\month\or January\or February\or March\or April\or
May\or June\or July\or August\or September\or October\or
November\or December\fi\space\number\day,\space\number\year}

\def\idag{\number\day. \ifcase\month\or januar\or februar\or marts\or
april\or maj\or juni\or juli\or august\or september\or oktober\or
november\or december\fi\space\number\year}

\newcount\minut
\newcount\our
\def\mintest{\ifnum\minut<10 {0\the\minut}\else \the\minut\fi}
\def\NU{\minut=\time \divide\time by60 \the\time.\multiply \time by60%
        \advance\minut by-\time \mintest}
\def\now{\minut=\time \divide\time by60 \our=\time%
        \ifnum\time<12  \the\time.\multiply \time by60%
        \advance\minut by-\time \mintest am%
        \else\ifnum\time>12 \advance\our by-12 \the\our.%
        \multiply \time by60 \advance\minut by-\time \mintest pm%
        \else\the\time.\multiply \time by60%
        \advance\minut by-\time \mintest pm%
        \fi\fi}

\edef\TheAtCode{\the\catcode`\@} \catcode`\@=11

\newcount\@aux@
\newcount\@month@
\newcount\@year@
\newcount\@wday@
\def\@weekday@{\@month@=\month \@year@=\year
        \ifnum\@month@=1 \@month@=13 \advance\@year@ by -1 \fi
        \ifnum\@month@=2 \@month@=14 \advance\@year@ by -1 \fi
        \@wday@=\day \advance \@wday@ by \@year@
        \@aux@=\@month@ \multiply\@aux@ by 2 \advance\@wday@ by \@aux@
        \@aux@=\@month@ \advance\@aux@ by 1 \multiply\@aux@ by 3
                \divide\@aux@ by 5
        \advance\@wday@ by \@aux@ \advance\@wday@ by 2
        \@aux@=\@year@ \divide\@aux@ by 4 \advance\@wday@ by \@aux@
        \@aux@=\@year@ \divide\@aux@ by 100 \advance\@wday@ by -\@aux@
        \@aux@=\@year@ \divide\@aux@ by 400 \advance\@wday@ by \@aux@
        \@aux@=\@wday@
        \divide\@aux@ by7 \multiply\@aux@ by7
    \advance\@wday@ by-\@aux@ }
\def\weekday{\@weekday@\ifcase\@wday@ Saturday\or Sunday\or Monday\or Tuesday\or
         Wednesday\or Thursday\or Friday\fi}
\def\ugedag{\@weekday@\ifcase\@wday@ L{\o}rdag \or S{\o}ndag \or
        Mandag \or Tirsdag \or Onsdag \or Torsdag \or Fredag \fi}

\def\NObreak{\penalty1000000 }

\newcount\ParanO
\newcount\ProcnO
\newcount\FignO
\def\advanceParan@O@{\global\advance\ParanO by 1
        \global\ProcnO=0 }
\def\advanceProcn@O@{\global\advance\ProcnO by 1 }
\def\advanceFign@O@{\global\advance\FignO by 1 }
\def\ProcnuM{\number\ProcnO}
\def\ParanuM{\number\ParanO}
\def\FignuM{\number\FignO}

\def\wh@{\eat@}
\def\wr@{\eat@}
\def\wa@{\eat@}

\newif\iffirstchapa@ \global\firstchapa@false


\def\chapclass{}
\outer\def\head#1\endhead{\gdef\chapclass{sec\-tion}%
  \penaltyandskip@{-200}\aboveheadskip
  \advanceParan@O@
  {\headfont@\raggedcenter@\interlinepenalty\@M
  \ParanuM. \,#1\endgraf}%
  \wh@{#1}
  \nobreak
  \vskip\belowheadskip}


\outer\def\appendix#1\endappendix#2{%
\gdef\chapclass{ap\-pen\-dix}%
  \penaltyandskip@{-200}\aboveheadskip
  \advanceParan@O@ \def\ParanuM{#2}
  {\headfont@\raggedcenter@\interlinepenalty\@M
  Appendix \ParanuM. \,#1\endgraf}%
  \wa@{#1}
  \nobreak
  \vskip\belowheadskip}

\outer\def\references#1\endreferences{%
\gdef\chapclass{the re\-fe\-ren\-ces}
  \penaltyandskip@{-200}\aboveheadskip
  {\headfont@\raggedcenter@\interlinepenalty\@M
  #1\endgraf}%
  \wr@{#1}
  \nobreak
  \vskip\belowheadskip}

\def\n{\advanceProcn@O@{\ParanuM.\ProcnuM}}
\def\tagn{\advanceProcn@O@\tag{\ParanuM.\ProcnuM}}
\def\figure{\advanceFign@O@{\smc Figure \FignuM. }}   

\def\output@{\shipout\vbox{%
 \iffirstpage@ \global\firstpage@false \global\firstchapa@false
  \pagebody \makefootline%
 \else \makeheadline \pagebody \makefootline
 \fi}%
 \advancepageno \ifnum\outputpenalty>-\@MM\else\dosupereject\fi}


\def\lmargnote#1{\setbox0=\hbox{\vbox{%
\hsize.7truein\raggedright\flushpar{\eightpoint#1}}}%
   \dimen0=\ht0\setbox1=\hbox{n}%
   \dimen1=\ht1\multiply\dimen1 by 2\advance\dimen0 by -\dimen1%
   \leavevmode\vadjust{\leftline{\smash{%
   \llap{\lower\dimen0\box0\quad}}}}\ignorespaces}
\def\rmargnote#1{\setbox0=\hbox{\vbox{%
\hsize.7truein\raggedright\flushpar{\eightpoint#1}}}%
   \dimen0=\ht0\setbox1=\hbox{n}%
   \dimen1=\ht1\multiply\dimen1 by 2\advance\dimen0 by -\dimen1%
   \leavevmode\vadjust{\rightline{\smash{%
   \rlap{\quad\lower\dimen0\box0}}}}\ignorespaces}

%
\def\ifundefined#1{\expandafter\ifx\csname#1\endcsname\relax}%
\def\ece#1#2{\expandafter#1\csname#2\endcsname}%

\newread\aux
\newwrite\aux

\newif\ifauxexist
\def\@readauxfile{%
   \if@auxfiledone \else
      \global\@auxfiledonetrue
      \immediate\openin\aux=\jobname.aux
      \ifeof\aux\immediate\closein\aux
         \message{\@undefinedmessage}%
         \global\@citewarningfalse
      \else\immediate\closein\aux
         \begingroup
            \endlinechar = -1
            \catcode`@ = 11
            \input \jobname.aux
         \endgroup
      \fi
      \immediate\openout\aux=\jobname.aux%
      \global\auxexisttrue%
   \fi
}%

\def\XREFS{
   \if@auxfiledone \else
      \global\@auxfiledonetrue
      \immediate\openin\aux=\jobname.aux
      \ifeof\aux\immediate\closein\aux
         \message{\@undefinedmessage}%
         \global\@citewarningfalse
      \else\immediate\closein\aux
         \begingroup
            \endlinechar = -1
            \catcode`@ = 11
            \input \jobname.aux
         \endgroup
      \fi
      \immediate\openout\aux=\jobname.aux%
      \global\auxexisttrue%
   \fi
}%
\newif\if@marglab
\global\@marglabtrue
\def\nomarglab{\global\@marglabfalse}

\newif\if@auxfiledone
\ifx\noauxfile\@undefined \else \@auxfiledonetrue\fi
\newwrite\aux
\def\@writeaux#1{\ifx\noauxfile\@undefined \write\aux{#1}\fi}%
\ifx\@undefinedmessage\@undefined
   \def\@undefinedmessage%
   {         Ingen .aux fil - jeg vil ikke advare om manglende labels.      }%
\fi

\newif\if@citewarning
\ifx\noauxfile\@undefined \@citewarningtrue\fi

\let\readauxfile = \@readauxfile
\let\writeaux = \@writeaux

\let\ifxrefwarning = \iftrue
\def\xrefwarningtrue{\@citewarningtrue \let\ifxrefwarning = \iftrue}%
\def\xrefwarningfalse{\@citewarningfalse \let\ifxrefwarning = \iffalse}%
\begingroup
  \catcode`\_ = 8
  \gdef\xrlabel#1#2{#1_x#2}%
\endgroup

\def\wrlabel#1#2#3{%
  \edef\temp{#1}%
  \@readauxfile
  \edef\@wl{\noexpand\writeaux{\string\@definelabel{\temp}{#2}{#3}}}%
  \@wl
  \ifdraft@\if@marglab\rmargnote{\temp}\fi\fi
  \ignorespaces
}%
\def\xwrlabel#1#2#3#4{%
  \edef\temp{#1}%
  \@readauxfile
  \edef\@wl{\noexpand\writeaux{\string\@definexlabel{\temp}{#2}{#3}{#4}}}%
  \@wl
  \ifdraft@\if@marglab\rmargnote{\temp}\fi\fi
  \ignorespaces
}%

\def\@definelabel#1#2#3{%
  \expandafter\gdef\csname\xrlabel{#1}{#2}\endcsname{#3}%
}%
\def\@definexlabel#1#2#3#4{%
  \expandafter\gdef\csname\xrlabel{#1}{#2}\endcsname{#4}%
  \expandafter\gdef\csname\xrlabel{#1}{#2cl}\endcsname{#3}
}%

\def\GeneralCheck#1#2{\readauxfile
  \expandafter \ifx\csname\xrlabel{#1}{#2}\endcsname\relax
    \if@citewarning
       \message{ Undefined label `#1' in line \the\inputlineno; }%
    \fi
    \expandafter\def\csname\xrlabel{#1}{#2}\endcsname{%
      `{\tt
        \escapechar = -1
        \expandafter\string\csname#1\endcsname
      }'%
    }%
  \fi
}%
\def\GeneralXCheck#1#2{\readauxfile
  \expandafter \ifx\csname\xrlabel{#1}{#2}\endcsname\relax
    \if@citewarning
       \message{ Undefined label `#1' in line \the\inputlineno; }%
    \fi
    \def\udskref##1##2{##2}
    \expandafter\def\csname\xrlabel{#1}{#2}\endcsname{%
      `{\tt
        \escapechar = -1
        \expandafter\string\csname#1\endcsname
      }'%
    }%
  \else\def\udskref##1##2{\leavevmode##1s~##2\ignorespaces}%
  \fi
}%

\def\norefwarnings{\global\@citewarningfalse}

\def\Pno{\ifnum\pageno<0\romannumeral-\pageno\else\the\pageno\fi}
\def\pagelab#1{\defineplab{#1}{p}{\noexpand\Pno}}%
\def\defineplab#1#2#3{%
  \edef\temp{#1}%
  \readauxfile
  \edef\@wr{\noexpand\writeaux{\string\@definelabel{\temp}{#2}{#3}}}%
  \@wr
  \ignorespaces
}%
\def\pagereftext#1{\def\pstring{#1}}
\pagereftext{page}                      
\def\pageref#1{\GeneralCheck{#1}p%
  \pstring~\csname\xrlabel{#1}p\endcsname%
}%
\def\pagerefs#1{\GeneralCheck{#1}p%
  \pstring s~\csname\xrlabel{#1}p\endcsname%
}%
\def\nopageref#1{\GeneralCheck{#1}p%
  \csname\xrlabel{#1}p\endcsname%
}%
\def\condpageref#1#2{\GeneralCheck{#1}p%
    \GeneralCheck{#2}p%
    \xdef\aref{\csname\xrlabel{#1}p\endcsname}%
    \xdef\bref{\csname\xrlabel{#2}p\endcsname}%
  \if\aref\bref%
    \pstring~\csname\xrlabel{#1}p\endcsname%
  \else%
    \pstring s~\csname\xrlabel{#1}p\endcsname\ and %
        \csname\xrlabel{#2}p\endcsname%
  \fi%
}%

\def\chapref#1{\if@citewarning
       \message{ Undefined label `#1' in line \the\inputlineno; }%
    \fi
    `{\tt #1}'}%
\def\chaprefs#1{\if@citewarning
       \message{ Undefined label `#1' in line \the\inputlineno; }%
    \fi
    `{\tt #1}'}%
\def\nochapref#1{\if@citewarning
       \message{ Undefined label `#1' in line \the\inputlineno; }%
    \fi
    `{\tt #1}'}%

\def\seclab#1{%
  \xwrlabel{#1}{s}{\chapclass}{\ParanuM}\ignorespaces}
\def\secref#1{\GeneralCheck{#1}s%
\leavevmode\csname\xrlabel{#1}{scl}\endcsname~\csname\xrlabel{#1}s\endcsname%
}%
\def\secrefs#1{\GeneralXCheck{#1}s
\udskref{\csname\xrlabel{#1}{scl}\endcsname}{\csname\xrlabel{#1}s\endcsname}%
}%
\def\nosecref#1{\GeneralCheck{#1}s%
    \csname\xrlabel{#1}s\endcsname%
}%

\def\numlab#1{%
  \wrlabel{#1}{n}{\ParanuM.\ProcnuM}\ignorespaces}
\def\numref#1{\GeneralCheck{#1}n%
  \csname\xrlabel{#1}n\endcsname%
}%
\def\tagref#1{\GeneralCheck{#1}n%
  \tagform@{\csname\xrlabel{#1}n\endcsname}%
}%

\def\proclab#1#2{%
  \xwrlabel{#1}{P}{#2}{\ParanuM.\ProcnuM}\ignorespaces}
\def\procref#1{\GeneralCheck{#1}P%
\leavevmode\csname\xrlabel{#1}{Pcl}\endcsname~\csname\xrlabel{#1}P\endcsname%
}%
\def\procrefs#1{\GeneralXCheck{#1}P%
\udskref{\csname\xrlabel{#1}{Pcl}\endcsname}{\csname\xrlabel{#1}P\endcsname}%
}%
\def\noprocref#1{\GeneralCheck{#1}P%
    \csname\xrlabel{#1}P\endcsname%
}%

\def\itemlab#1{%
  \edef\temp{#1}%
  \@readauxfile
  \edef\@wl{\noexpand\writeaux%
  {\string\@definelabel{\temp}{i}{\the\rostercount@}}}%
  \@wl
  \ignorespaces%
}
\def\marglab#1{\ifdraft@\if@marglab\rmargnote{#1}\fi\fi\ignorespaces}
\def\itemreftext#1{\def\istring{#1}}
\itemreftext{item}                              
\def\itemref#1{\GeneralCheck{#1}i%
  \istring~\rom{(\csname\xrlabel{#1}i\endcsname)}%
}%
\def\itemrefs#1{\GeneralCheck{#1}i%
  \istring i~\rom{(\csname\xrlabel{#1}i\endcsname)}%
}%
\def\noitemref#1{\GeneralCheck{#1}i%
  \rom{(\csname\xrlabel{#1}i\endcsname)}%
}%


\def\figreftext#1{\def\fstring{#1}}
\figreftext{figure}
\def\figref#1{\GeneralCheck{#1}f%
  \fstring~\csname\xrlabel{#1}f\endcsname%
}%
\def\nofigref#1{\GeneralCheck{#1}f%
  \csname\xrlabel{#1}f\endcsname%
}%


\def\divref#1{\GeneralCheck{#1}d%
  \ignorespaces\csname\xrlabel{#1}d\endcsname%
}%
\def\xdivref#1{\GeneralXCheck{#1}d%
\leavevmode\unskip\csname\xrlabel{#1}{dcl}\endcsname~\csname%
\xrlabel{#1}d\endcsname%
}%

\newread\nix
\newwrite\nix

\newif\ifnixexist
\def\@readnixfile{%
   \if@nixfiledone \else
      \global\@nixfiledonetrue
      \immediate\openin\nix=\jobname.nix
      \ifeof\nix\immediate\closein\nix
         \message{\@undefinedmessagenix}%
         \global\@citenixwarningfalse
      \else\immediate\closein\nix
      \fi
      \immediate\openout\nix=\jobname.nix%
      \global\nixexisttrue%
   \fi
}%

\newif\if@nixfiledone
\ifx\nonixfile\@undefined \else \@nixfiledonetrue\fi
\newwrite\nix
\def\@writenix#1{\ifx\nonixfile\@undefined \write\nix{#1}\fi}%
\ifx\@undefinedmessagenix\@undefined
   \def\@undefinedmessagenix%
   {         Ingen .nix fil - jeg vil ikke advare om manglende indeceringer.      }%
\fi

\newif\if@citenixwarning
\ifx\nonixfile\@undefined \@citenixwarningtrue\fi

\let\readnixfile = \@readnixfile
\let\writenix = \@writenix

\let\ifxrefnixwarning = \iftrue
\def\xrefnixwarningtrue{\@citenixwarningtrue \let\ifxrefnixwarning = \iftrue}%
\def\xrefnixwarningfalse{\@citenixwarningfalse \let\ifxrefnixwarning = \iffalse}%

\def\wrnotidx#1#2{%
  \edef\temp{#1}%
  \@readnixfile
  \edef\@wl{\noexpand\writenix{\string\@definenotidx{\temp}{#2}}}%
  \@wl%
  \ignorespaces
}%

\def\leaderfill{\unskip\nobreak\penalty50
        \hskip0pt\hbox{}\nobreak\leaders\hbox to 5pt{\hss.\hss}\hfill}

\def\@definenotidx#1#2{%
  \csname\xrlabel{#1}i\endcsname \leaderfill #2\cr%
}%

\def\norefnixwarnings{\global\@citenixwarningfalse}

\def\GeneralnixCheck#1#2{\readauxfile
  \expandafter \ifx\csname\xrlabel{#1}{#2}\endcsname\relax
    \expandafter\def\csname\xrlabel{#1}{#2}\endcsname{{\bf ?}}%
  \fi
}%
\def\nixref#1{\GeneralnixCheck{#1}p%
  \csname\xrlabel{#1}p\endcsname%
}%

\newcount\ntidx \ntidx=1
\def\numlet#1{\romannumeral#1}

\def\notidx#1#2{\definenix{#1}{#2}{\noexpand\Pno}}%
\def\definenix#1#2#3{%
  \edef\temp{idx_\numlet{\ntidx}}%
  \readnixfile
  \edef\@wr{\noexpand\writenix{\string\@definenotidx{\temp}{#3}}}%
  \@wr
  \expandafter\gdef\csname\xrlabel{idx_\numlet{\ntidx}}i\endcsname{#1&#2}%
  \advance\ntidx by 1%
  \ignorespaces
}%

\def\xnotidx#1#2#3{\definexnix{#1}{#2}{\noexpand\Pno}{#3}}%
\def\definexnix#1#2#3#4{%
  \edef\temp{idx_\numlet{\ntidx}}%
  \def\xPno{\noexpand#4}%
  \readnixfile
  \edef\@wr{\noexpand\writenix{\string\@definenotidx{\temp}{#3, \xPno}}}%
  \@wr
  \expandafter\gdef\csname\xrlabel{idx_\numlet{\ntidx}}i\endcsname{#1&#2}%
  \advance\ntidx by 1%
  \ignorespaces
}%

\let\notidxfont@=\refsfont@
\newdimen\colone%
\newdimen\coltwo%
\def\widestentry#1{\setbox0=\hbox{\notidxfont@ #1\quad}\dimen0=\wd0%
        \colone=0pt \advance\colone by\dimen0%
        \coltwo=6truein \advance\coltwo by -\colone }

\outer\def\notaindex#1\endnotaindex{%
\gdef\chapclass{}
  \penaltyandskip@{-200}\aboveheadskip
  {\headfont@\raggedcenter@\interlinepenalty\@M
  #1\endgraf}%
  \wr@{#1}
  \nobreak
  \vskip\belowheadskip
  \immediate\closeout\nix
  \advance\ntidx by -1
  \immediate\write16{^^J Notations indgange skrevet: \the\ntidx ^^J}
         \begingroup \notidxfont@
            \endlinechar = -1
            \catcode`@ = 11
            \halign{\hbox to \colone{##\hfil}&%
                \hbox to \coltwo{##}\cr
            \input \jobname.nix
            }
         \endgroup
  }

\colone=1truein \coltwo=5truein

%
\def\leaderfill{\unskip\nobreak\penalty50
        \hskip0pt\hbox{}\nobreak\leaders\hbox to 5pt{\hss.\hss}\hfill}

\def\hpost#1#2#3{#1\ignorespaces#2\ignorespaces
\leaderfill#3\ignorespaces}
\def\Post#1#2#3#4{\ifcase#1 \or\hpost#2#3#4\else\fi}
\def\Read{\ifeof\toc\immediate\closein\toc
        \else\read\toc to\tocslags\ifeof\toc\immediate\closein\toc
                \else\read\toc to\tocnr
                     \read\toc to\tocnavn
                     \read\toc to\toctal
                     \Post{\tocslags}{\tocnr}{\tocnavn}{\toctal}\par
                     \Read\fi\fi}

\newread\toc
\newwrite\toc
\newif\ifcontent@ 
\outer\def\contents#1\endcontents{%
        \global\content@true 
        \gdef\wh@##1{{\let\the=0\let\romannumeral=0\edef\next{%
                \write\toc{1}%
                \write\toc{{\ParanuM.\enspace}}%
                \write\toc{{##1}}%
                \ifnum\pageno<0\write\toc{{\romannumeral-\pageno}}
                        \else\write\toc{{\the\pageno}}\fi}\next}}
        \gdef\wa@##1{{\let\the=0\let\romannumeral=0\edef\next{%
                \write\toc{1}%
                \write\toc{{Appendix \ParanuM.\enspace}}%
                \write\toc{{##1}}%
                \ifnum\pageno<0\write\toc{{\romannumeral-\pageno}}
                        \else\write\toc{{\the\pageno}}\fi}\next}}
        \gdef\wr@##1{{\let\the=0\let\romannumeral=0\edef\next{%
                \write\toc{1}%
                \write\toc{{}}%
                \write\toc{{##1}}%
                \ifnum\pageno<0\write\toc{{\romannumeral-\pageno}}
                        \else\write\toc{{\the\pageno}}\fi}\next}}
        \penaltyandskip@{-200}\aboveheadskip
          {\headfont@\raggedcenter@\interlinepenalty\@M
          #1\endgraf}%
          \nobreak
          \vskip\belowheadskip
        {\parindent=0pt%
        \immediate\openin\toc=\jobname.toc
        \Read\par}
        \immediate\openout\toc=\jobname.toc
        \bigskip\bigskip}

\font\headf=cmcsc8 \footline{\hfil}
\newif\ifdraft@
\def\draft{\global\draft@true%
        \headline{\headf \hss Draft, \weekday\ \today, \now\hss}
        \footline{\hss\tenrm\folio\hss}
        \def\logo@{\baselineskip2pc \hbox to\hsize{\hfil\eightpoint\smc
        Draft, \weekday\ \today, \now}}}
\def\nohead{\headline{\hfil}
        \footline{\hss\tenrm\the\pageno\hss}}

\ParanO=0 \FignO=0

\outer\def\enddocument{%
\ifmonograph@ 
\else
 \nobreak
 \thetranslator@
 \count@\z@ \loop\ifnum\count@<\addresscount@\advance\count@\@ne
 \csname address\number\count@\endcsname
 \csname email\number\count@\endcsname
 \repeat
\fi
 \vfill\supereject%
 \ifcontent@%
      \ifauxexist\immediate\write16{^^J^^J^^J
 HUSK at TeX'e to gange, da referenser og indholdsfortegnelsen kan have aendret sig!!!
 ^^J^^J^^J}\immediate\closeout\toc\immediate\closeout\aux%
      \else\immediate\write16{^^J^^J^^J
 HUSK at TeX'e to gange, da indholdsfortegnelsen kan have aendret sig!!!
 ^^J^^J^^J}\immediate\closeout\toc\fi%
 \else\ifauxexist\immediate\write16{^^J^^J^^J
 HUSK at TeX'e to gange, da krydsreferenserne kan have aendret sig!!!
 ^^J^^J^^J}\immediate\closeout\aux%
      \fi%
 \fi\end}

\catcode`\@=\TheAtCode\relax

\def\proof{\demo{Proof}}
\define\endproof{\hfill\penalty50
    \hbox{}\nobreak\hfill$\square$\enddemo}

\input xy
\xyoption{arrow}\xyoption{matrix}\xyoption{tips}\xyoption{curve}

\SelectTips{cm}{}


\input texdraw
\XREFS
\nomarglab

\loadeusm

\define\Phf{\eusm P}

\define\llangle{\langle\kern-2.5pt\langle}
\define\rrangle{\rangle\kern-2.5pt\rangle}

\define\lra{\longrightarrow}
\define\ra{\to}
\define\Emu{\Cal E_\mu}
\define\E{\Cal E}

\define\C{\Bbb C}
\define\Z{\Bbb Z}
\define\M{M}
\define\bM{\bar \M}

\define\bDh{\bar \Delta^{1/2}}
\define\Dc{\Cal D}
\define\LD{{\Lambda}_{\Delta}}

\define\F{\Cal F}
\define\tF{\tilde {\F}}
\define\Pa{P_a}
\define\tPa{\tilde P_a}
\define\tP{\tilde {\Cal P}}
\define\Xt{\check{X}}
\define\xt{\check{x}}
\define\Lt{\check{L}}
\define\bH{\bar H}
\define\pit{\check{\pi}}
\define\Phit{\check{\Phi}}
\define\aut{\text{Aut}}
\define\ttau{\tilde\tau}

\define\conge{\cong_e}

\define\sign{\text{sign}}

\define\Li{\Cal L_\Dc}

\define\res{\operatorname{Res}}

\define\Mll#1#2{M^{#1}_{#2}}
\define\Ml#1{\Mll{#1}{\langle L\rangle}(Y)}

\define\pr{\operatorname{pr}}
\define\Gr{\operatorname{Gr}}
\define\End{\operatorname{End}}
\define\Aut{\operatorname{Aut}}
\define\Hom{\operatorname{Hom}}
\define\GL{\operatorname{GL}}
\define\PGL{\operatorname{PGL}}
\define\Pj{\Bbb P}
\define\lcm{\operatorname{lcm}}

\define\Pic{\operatorname{Pic}}
\define\Div{\operatorname{Div}}
\define\ord{\operatorname{ord}}
\define\Nm{\operatorname{Nm}}
\define\Tr{\operatorname{Tr}}
\define\Trt#1{\Tr\big(\tau \big| H^{#1}\big(M(X),\Cal L^\kappa\big)\big)}
\define\Trtk#1{\Tr\big(\tau \big| H^{#1}\big(M(X;k),\varphi^*\Cal L^\kappa\big)\big)}
\define\id{\operatorname{Id}}
\define\pa{\operatorname{par}}
\define\so{\Cal O}

\define\en#1{e^{\frac{2\pi i}n #1}}

\define\em{\Cal F}

\define\elm{\Gamma}

\define\deq{=}

\define\e#1{\rlap{$\pmb #1$}\hphantom{#1}}
\define\sse#1{\rlap{$\ssize\pmb#1$}\hphantom{\ssize#1}}

\define\Spec{\operatorname{Spec}}

\define\im{\operatorname{Im}}

\loadeusm
\define\SExt{\eusm E\kern-0.6ptxt}
\define\SHom{\eusm H\kern-0.6ptom}

 \mathchardef\mycall="324C 
 \mathchardef\mycale="3245 
 \mathchardef\mycalm="324D 
 \mathchardef\mylbrace="3266
 \mathchardef\myrbrace="3267

 \newcount\refno
 \refno=0
 \define\defref#1.{\advance\refno by 1 %
   \global\edef#1{\the\refno}}

 \defref\refAn.
 \defref\refAU.
 \defref\refAtone.
 \defref\refAxDW.
 \defref\refBFM.
 \defref\refBFQ.
 \defref\refBL.
 \defref\refF.
 \defref\refG.
 \defref\refGrfour.
 \defref\refGr.
 \defref\refH.
 \defref\refJKKW.
 \defref\refKo.
 \defref\refL.
 \defref\refLa.
 \defref\refMeS.
 \defref\refMuone.
 \defref\refMutwo.
 \defref\refNR.
 \defref\refQ.
 \defref\refSe.
 \defref\refTh.
 \defref\refTUY.
 \defref\refW.
 \defref\refWi.


 \topmatter

\title Automorphism fixed Points in the Moduli Space of Semistable bundles\endtitle

 \author J{\o}rgen Ellegaard Andersen and Jakob Grove\endauthor

 \thanks This research was partially conducted
 for the Clay Mathematics Institute at University of California, Berkeley by the first author.
 The second author was supported by Japan Society for the Promotion of Science.\endthanks

 \address J.E.~Andersen, Mathematics Department, University of Aarhus, DK-8000 Aarhus C,
 Denmark \endaddress

 \curraddr University of California, Berkeley, Berkeley, CA 94720,
 USA\endcurraddr

\email andersen\@imf.au.dk\endemail

\address J.~Grove, Department of Mathematics Graduate School of Science, Kyoto
 University, Ki\-ta\-shi\-ra\-ka\-wa, Sakyo-ku, Kyoto 606-8502, Japan     \endaddress

 \email grove\@kusm.kyoto-u.ac.jp, grove\@imf.au.dk\endemail


 \abstract Given an automorphism $\tau$ of a smooth complex algebraic curve $X$, there is an
 induced action on the moduli space $\M$ of semi-stable rank $2$
 holomorphic bundles with fixed determinant.
 We give a complete description of the fixed variety in terms of moduli
 spaces of parabolic bundles on the quotient curve $X/\langle \tau \rangle$.\endabstract

 \keywords Fixed points, Moduli spaces of holomorphic bundles, Gauge theory approach to
 quantum invariants of $3$-manifolds\endkeywords

 \subjclass 14D20, 
    14H60, 
    14H37, 
    32G13, 
    57M27  
 \endsubjclass

 \endtopmatter

\document
 \contents Contents\endcontents


\head Introduction\endhead

Throughout this paper $X$ will be a smooth complex algebraic
curve, and $\tau$ and automorphism of $X$, say of order $n$. Let
$Y = X / \langle\tau \rangle$ be the quotient curve, which is of
course again a smooth complex algebraic curve. Let $\pi : X \ra Y$
be the corresponding possibly ramified covering. We will assume
that $Y$ is connected. However, any statement may be easily
generalised by considering one component of $Y$ at a time.

Let $M$ be the moduli space of semi-stable holomorphic bundles of
rank $2$ and fixed determinant on $X$. We get an induced action of
$\langle\tau \rangle$ on $M$, provided that this fixed determinant
line bundle is preserved up to isomorphism under pullback by
$\tau$. Let us denote the fixed variety of this group action by
$|M|$.

In the present paper we give a construction of the fixed variety
$|M|$ from a certain finite set of admissible moduli spaces of
semi-stable parabolic bundles on $Y$ with the parabolic structures
concentrated at the ramification points of $\pi \: X \ra Y$.

The construction in outline is as follows: First we establish that
any fixed point in the moduli space $M$, can be represented by an
equivariant (semi-stable) bundle, i.e. by a bundle with an action
of the group $\langle\tau \rangle$ covering the action on $X$. By
fixing the determinant line bundle equivariantly, to the extend
this can be done, we reduce the ambiguity in the choice of the
corresponding equivariant bundle as much as possible (see
\procref{repfixeb} and the paragraphs following it).

This prompts us to study equivariant bundles in general up to
equivariant isomorphism. Suppose $(E,\ttau)$ is an equivariant
bundle. In the case where $\pi: X \ra Y$ is ramified, there are
obvious numeric invariants associated to the equivariant bundle
$(E,\ttau)$ obtained in the following manner: At a special orbit
$y$ of $\tau$, say of length $k(y)$, $\ttau^{k(y)}$ acts on the
fiber $E_x$ over each point $x$ in the orbit. Thus, its
eigenvalues, which are $\frac n{k(y)}$-roots of unity, will be
invariants of the equivariant isomorphism class of $(E,\ttau)$.
These form discrete invariants of the equivariant bundle. To
determine the equivariant bundle up to equivariant isomorphism, we
devise a scheme for performing a number of equivariant elementary
modifications at the special orbits, so as to obtain a
quasi-parabolic bundle, which is seen to be equivariantly
isomorphic to the pullback of a quasi-parabolic bundle $(\bar E,
\bar F)$ on $Y$. We specify weights $w$ for this quasi-parabolic
bundle $(\bar E, \bar F)$ on $Y$ as functions of the discrete
eigenvalue invariants of $(E, \ttau)$.

The next step is to specify a set $\tilde P_a$ of admissible
parabolic bundles (see \procref{labelset} and
\noprocref{admisparabb}), and prove that these are exactly the
parabolic bundles obtained as describe above. This is done in
\procref{bijection}, where we establish that there is a bijective
correspondence between this set of admissible parabolic bundles
and then the set of equivariant bundles as specified above. We
view this bijection as a construction of all such equivariant
bundles.

This construction behaves very well with respect to
semi-stability. Namely, an admissible parabolic bundle is
parabolic semi-stable if and only if the corresponding equivariant
bundle is semi-stable as a vector bundle as stated in
\procref{ssiffpss}. Furthermore, S-equivalent parabolic bundles
are taken to S-equivalent bundles, hence we get a well defined set
map $\F$ from the moduli space of admissible semi-stable parabolic
bundles $P_a$ on $Y$ to the fixed variety $|M|$. --- Since any
fixed point can be represented by an equivariant bundle, this map
is clearly surjective. The way we have defined the set of
admissible parabolic bundles, implies that this map has finite
fibers and we give a complete description of all fibers of $\F$ in
the main set-theoretic \procref{mainsetT}. The content of this
theorem is in words:

There is a specific equivalence relation on admissible
semi-stable, but not stable, parabolic bundles, induced from a
finite group action on parabolic line bundles, which describes the
fibers of $\F$. If $n$ is odd, $\F$ is a bijection between stable
admissible parabolic bundles and stable fixed points and it
induces a bijection between equivalence classes of semi-stable
non-stable parabolic bundles and semi-stable non-stable fixed
points.

If $n$ is even, there is an involution on the moduli space of
admissible parabolic bundles under which $\F$ is invariant. If the
greatest common divisor $r$ of the orbit lengths for $\tau$ is
odd, $\F$ gives a bijection between the quotient of the moduli
space of stable admissible parabolic bundles and the stable fixed
points, and again there is a bijection between equivalence classes
of semi-stable non-stable parabolic bundles and semi-stable
non-stable fixed points.

If $r$ is even, $\F$ induces a bijection between the quotient of
the stable admissible parabolic bundles which are not fixed by the
involution and stable fixed points. The stable admissible
parabolic bundles, fixed by the involution together with the
equivalence classes of the semi-stable but not stable admissible
parabolic bundles, goes bijectively under $\F$ into the
semi-stable but non-stable points in the fixed variety.

We describe in details how these equivalence relations affect the
discrete invariants, i.e.\ the parabolic weight and how they are
induced by taking elementary modifications and tensoring with
certain line bundles (see point 1.\ and 2.\ following
\procref{mainsetT}).

Using the geometric invariant theory construction of these moduli
spaces, we analyse our set-theoretic map $\F$ and show that it is
a morphism of varieties in \procref{mod sp morphism}. Furthermore,
\procref{em normal even} establishes that $\F$ (or in the case of
$n$ even, the induced map on the $\Bbb Z/2$-quotient mentioned
before) is a birational equivalence, which, when restricted to
components of the moduli space of admissible parabolic bundles
(resp.\ its quotient), is the normalising map onto the
corresponding irreducible components of the fixed variety.

As examples, the un-ramified case and ramified hyper-elliptic case
is treated in details.

The motivation for an explicit construction of this fixed variety
in terms of parabolic bundles stems from the gauge theoretic
approach to $2+1$-dimensional topological quantum field theories,
as initially outlined by Witten in \cite{\refWi}, Axelrod, Della
Pietra and Witten in \cite{\refAxDW} and Atiyah in
\cite{\refAtone}. From the algebraic geometric viewpoint, this
approach has received considerable attention, e.g. Hitchin
\cite{\refH}, Faltings \cite{\refF} and Thaddeus \cite{\refTh},
Beauville \& Laszlo \cite{\refBL}, Narasimhan and Ramadas
\cite{\refNR} to mention a few. The program of proving all the
axioms of a full TQFT, has however not been complete from this
purely gauge-theoretic point of view. On the conformal field
theory side, further progress was made in the direction of
providing a full construction of a modular functor, which in turn
will provide the basis for another combinatorial construction of
these TQFTs, namely the construction of TQFT from Modular
functors, see \cite{\refW} and \cite{\refG}. Here the works of
Segal \cite{\refSe} and Tsuchiya, Ueno and Yamada \cite{\refTUY}
stand out. Work is in progress to provide a full construction of a
modular functor based solely on the techniques of \cite{\refTUY}
(see \cite{\refAU}). Once this is complete this can be combined
with the results of Laszlo \cite{\refLa}, to establish that the
association of the vector space $H^0(M,\Theta)$, where $\Theta \in
\Pic(M)$ is ample, to the curve $X$, can be extended to a full
modular functor.

From this point of view it is therefore an interesting problem to
use the geometry of the moduli space $M$ to study the character of
the resulting representation of the automorphism group of $X$. As
explained in \cite{\refAn}, this character can be expressed as a
certain cohomological/homological pairing on the fixed variety
$|M|$, by an application of the Lefschetz-Riemann Rock formula due
to Baum, Fulton, Quart and MacPherson, \cite{\refBFM},
\cite{\refBFQ} and \cite{\refQ}. This gave a proof of the
Asymptotic Expansion Conjecture for mapping tori of automorphisms
of curves. --- To further study this character, a complete
description of the fixed variety is needed. Such a description is
given in this paper as explained above in terms of the moduli
space of admissible parabolic bundles $P_a$. It is therefore an
interesting problem to pull back this pairing to the moduli space
of admissible parabolic bundles bundles $P_a$ on $Y$ and seek an
evaluation there of the pairing, e.g. by means of the results of
Jeffrey, Kiem, Kirwan and Woolf, \cite{\refJKKW}, on the
intersection cohomology of these moduli spaces.

We remark, that although the techniques presented here all
generalise to higher rank, only rank $2$ is treated. This choice
was made in order to eliminate the more involved bookkeeping
required to treat the general rank case.

Both authors wish to express their appreciativeness for the
hospitality of the Mittag-Leffler Institute and the Department of
Mathematics at Kyoto University, where a significant part of the
work at hand was carried out. We also which to thank K.~Ueno for
helpful discussions on this project.

\head Equivariant Bundles\endhead

\seclab{prelim}

Recall that $X$ is a smooth complex algebraic curve, $\tau \: X
\to X$ an automorphism of order $n$, which gives a possibly
ramified covering $\pi \: X \to Y = X/\langle\tau\rangle$ of
curves.

In the sequel we employ the convention that the set of special
orbits for $\tau$ is denoted $P\subset Y$. For every orbit $y \in
P$ we denote its length by $k = k(y)$, so the ramification number
is $n'=n'(y)= \frac n{k(y)}$. Let $\xi = \xi_y$ be the $n'$-th
root of unity given by $\xi_y= d_x \tau^k: T_x X \ra T_x X$ for
any $x\in \pi^{-1}(y)$. It is not hard to see that there exists a
neighbourhood $U = U_x$ of $x$ in which there are local
coordinates $z$ centered around $x$ such that $\tau^k$ can be
expressed as $z \circ \tau^{-k} = \xi^{-1}_y \cdot z$.

Before we go on, observe the following important remark:

\remark{Remark \n\proclab{cover factor}{Remark}} For any divisor
$d$ of $n$, the length of the orbit of $\tau^d$ through a point $x
\in X$ is $k\big/\gcd\big(d,k\big)$. This is because the length of
the orbit is the smallest $l > 0$ such that $k$ divides $dl$,
which is the factor of $k$ that is not in $d$:
$\lcm\big(d,k\big)\big/d = k\big/\gcd\big(d,k\big)$. This means
that the length of the fiber $\pi_d^{-1}\big(\pi(x)\big)$ of the
projection $\pi_d \: X/\langle\tau^d\rangle \to Y$ is
$\gcd\big(d,k\big)$.

The projection $\pi_d \: X/\langle\tau^d\rangle \to Y$ is
unramified precisely when the lengths of its fibers are constantly
equal to the generic length $d$. That in turn means that
$\gcd\big(d,k\big) = d$ for all $x$, in particular, that $d$
divides all $k$, i.e.\ that $d$ divides the greatest common
divisor $r = \gcd \{\, k(y) \mid y\in Y\, \} $ of the orbit
lengths.

This shows that any ramified cover coming from the action of a
single automorphism $\tau$ of $X$ factors into a ramified part
$\pi_r \: X \to \hat X = X/ \langle\tau^r\rangle$ for which the
fiber lengths are co-prime, and an unramified part $\pi_u \: \hat
X \to Y$ whose degree is the greatest common divisor $r$ of the
orbits lengths of $\tau$.\endremark

Let now $\M$ be the moduli space of S-equivalence classes of
semi-stable bundles on $X$ of rank $2$ and determinant isomorphic
to some fixed line bundle. Provided that this line bundle is
invariant under pullback by $\tau$, we see that pullback by $\tau$
induces an action of $\langle\tau\rangle$ on the set of
S-equivalence classes of such bundles. It is easily seen that this
action is indeed well defined, for if $W_i$, $i=1,2$, are two
semi-stable holomorphic bundles, such that $\Gr(W_1)\cong
\Gr(W_2)$ then
 $$
   \Gr(\tau^* W_1) \cong \tau^*\Gr(W_1) \cong
        \tau^*\Gr(W_2) \cong \Gr(\tau^*W_2).
 $$
As pullback by morphisms between smooth complex algebraic curves
gives morphisms between the moduli spaces, we actually get an
induced algebraic action of $\langle\tau\rangle$ on $\M$ (see
section 4) and we denote the fixed variety by $|\M|$. Note that a
bundle $W$ represents a point in $|\M|$ if and only if
$\Gr\big(\tau^*(W)\big) \cong \Gr(W)$. However, this same point in
$|\M|$ is also represented by $\Gr(W)$ which satisfies that
$\tau^*\big(\Gr(W)\big) \cong \Gr(W)$. Hence, any fixed point in
$\M$ can be represented by a semi-stable bundle which is preserved
under pullback by $\tau$.

Consider now a holomorphic vector bundle $V$ over $X$ of any rank,
and suppose that $\tau^*V$ is isomorphic to $V$. Then there exists
a bundle map $\tilde \tau : V \to V$ covering $\tau$. Two such
pairs $(V_\nu,\tilde\tau_\nu)$ are said to be isomorphic, if there
is an isomorphism of bundles $\Phi : V_1 \to V_2$ intertwining
$\tilde{\tau}_1$ and $\tilde{\tau}_2$ and we write $(V_1,\ttau_1)
\cong (V_2,\ttau_2)$ or just $V_1 \conge V_2$.

\proclaim{Definition \n} A pair $(V,\tilde\tau)$ consisting of a
holomorphic vector bundle $V$ and a bundle map $\tilde\tau$
covering $\tau$ is called a lift, if $$ \tilde\tau^n = \id_V. $$
\endproclaim

We shall also refer to such a pair $(V,\tilde\tau)$ as an {\it
equivariant holomorphic bundle}\/. We note that the group of
$n$-th roots of unity $\zeta_n$ acts freely on the set of
isomorphism classes of equivariant holomorphic bundles.

If we are given a divisor on $X$ which is $\tau$-invariant, then
there is a naturally induced lift to the bundle associated to the
divisor, just by inducing the action from the
$\langle\tau\rangle$-action on the meromorphic functions ${\Cal
M}(X)$ on $X$. Here we use the convention for $\tau^* : {\Cal
M}(X) \ra {\Cal M}(X)$ that $\tau^* f = f \circ \tau^{-1}$.

If $V$ is a simple bundle which is preserved by $\tau$, i.e.
$\tau^*V \cong V$, then we can always find a lift of $\tau$ to $V$
and $\zeta_n$ acts simply transitive on the set of such lifts. In
particular, this is the case for all preserved stable bundles and
in particular for all preserved line bundle.\footnote{We shall
actually see shortly that any preserved line bundle can be
represented by a $\tau$-invariant divisor.}

This may not always be the case for preserved semi-stable vector
bundles. However, if we only consider them modulo S-equivalence,
then, as we shall see now, this can always be achieved: Assume
that a semi-stable bundle $W$ of rank $2$ on $X$ represents a
fixed point in $|\M|$, i.e. it satisfies that its S-equivalence
class is preserved by $\tau$
 $$
   \Gr(\tau^*W)  \cong \Gr(W).
 $$
It is in this case easily shown that there is a lift $\tilde\tau$
of $\tau$ to $\Gr(W)$:

The graded object of $W$ is $\Gr(W) = L_1 \oplus L_2$, where
$L_\nu$ is a line bundle of degree $\frac12\deg E$ for $\nu = 1,
2$. So $\tau^*(L_1 \oplus L_2) \cong L_1 \oplus L_2$ and since any
non-zero homomorphism between line bundles of the same degree is
an isomorphism, we either get that $\tau^*L_\nu \cong L_\nu$
(which we refer to as the {\it invariant\/} case) or $\tau^*L_1
\cong L_2$ and vice versa, but $L_1 \ncong L_2$ (which we refer to
as the {\it degenerate\/} case). In the later case we observe that
$(\tau^*)^2 L_\nu \cong L_\nu$, which in particular means that $n$
must be even, for this case to occur.

In the invariant case, it follows from the above argument, that we
may choose $(L_\nu,\tilde\tau_\nu)$ so that $\tilde\tau_\nu^n =
\id_{L_\nu}$. For the degenerate case suppose $\tilde\tau :
L_1\oplus L_2 \ra L_1\oplus L_2$ is a bundle isomorphism covering
$\tau$. Then, since $L_1 \ncong L_2$, $\ttau$ must be
off-diagonal. However $\tilde \tau^n$ has to be diagonal with
respect to the splitting, since $H^0\big(\End(L_1\oplus L_2)\big)
= \Bbb C \oplus \Bbb C$, and a moments thought shows one can find
a diagonal matrix $\Lambda$ so that $(\Lambda \tilde \tau)^n =
\id_{L_1 \oplus L_2}$.

Hence, we have now seen that any fixed point in $\M$ can be
represented by an equivariant bundle $(E,\tilde\tau)$.

Now observe that if we change a lift $\tilde\tau$ by a $\mu\in
\zeta_n$, then the induced lift to the determinant $\det(E)$ is
changes by $\mu^2$. We can therefore always change the lift
$\tilde \tau$ so that the induced lift on the determinant is a
given one in the case $n$ is odd and one of two possibilities in
case $n$ is even. Fix therefore a set $\Dc$ of such equivariant
line bundles whose underlying bundles is the determinant fixed
above. I.e.\ in the case $n$ is odd, $\Dc$ contains exactly one
lift to the determinant line bundle and in case $n$ is even
exactly two lifts to the determinant line bundle, which are not
equivalent under the action of $\{\, \mu^2\mid \mu\in
\zeta_n\,\}$. We have thus arrived at

\proclaim{Lemma \n\proclab{repfixeb}{Lemma}} Every fixed point in
$\M$ can be represented by a semi-stable equivariant bundle
$(E,\tilde\tau)$ with $\det(E,\tilde\tau)\in \Dc$. \hfill$\square$
\endproclaim

We observe, in the case of stable bundles, that when $n$ is odd,
we have by these means a unique way to represent such fixed points
by equivariant bundles. When $n$ is even there are exactly two
such equivariant bundles, say $(E,\tilde\tau)$ and then
$(E,-\tilde\tau)$. For semi-stable, but not stable bundles, the
situation is of course more involved. --- We shall see the
significance of fixing the determinant equivariantly in the
following section.

Let $\Li$ be the set of isomorphism classes of equivariant bundles
$(E,\tilde\tau)$ with $\det(E,\tilde\tau) \in \Dc$, $\Li'$ the
subset consisting of isomorphism classes of equivariant bundles,
where the underlying bundle is semi-stable and let $\Pi_e : \Li'
\ra |\M|$ be the projection map, which forgets the lift and takes
the S-equivalence class of the underlying bundle. The content of
\procref{repfixeb} is exactly that $\Pi_e$ is surjective.

Let now $(E,\ttau)$ be an equivariant bundle. For each $y\in P$ we
now define two integers $0\leq d_1=d_1(y) \leq d_2=d_2(y)
<n'$ by requiring that
 $$
   \theta_\nu = \xi^{d_\nu}, \qquad\nu =1,2,
 $$
where $\theta_\nu = \theta_\nu(y)$, $\nu=1,2$, are the eigenvalues
of $\tilde\tau^k$ acting on fibers of $E$ over $\pi^{-1}(y)$. We
note that the ordered pair $(d_1,d_2)$ is an invariant of the
isomorphism class of the equivariant bundle $(E,\ttau)$. We call
this pair the {\it numeric data\/} or {\it numeric invariant\/} of
$(E,\ttau)$, and it is a map $$
 (d_1,d_2): \Li \lra \prod_{y\in P}T_{n'(y)},
$$ where
 $$
  T_n =\{(d_1,d_2)\in \Z \mid 0 \leq d_1 \le d_2 < n\}.
$$ The numeric data is a discrete invariant of an equivariant
bundle, but we need of course much more than this discrete
invariant to determine the equivariant bundle. We proceed as
follows in our further analysis of $(E,\ttau)$.

Let $\E$ be the sheaf of sections in $E$. Let $\mu \in \zeta_n$
and define the eigensubsheaf $\Emu$ of $\pi_*\E$ corresponding to
the eigenvalue $\mu$ by
 $$
   \Emu(U) = \big\{\,s\in \E\big(\pi^{-1}(U)\big) \,\big|\,
   \tilde \tau s = \mu s\,\big\}
 $$
for any open subset U of $Y$, where $(\ttau s)(x) =
\ttau\big(s(\tau^{-1}(x))\big)$, $x\in X$.

 Let us now consider the sheaf
morphisms $\pi_\mu \: \pi_*\E \to \pi_*\E$ by letting
 $$
   \pi_\mu = \frac1n \sum_{l=0}^{n-1}
   \mu^{-l}\tilde\tau^l. \tagn
 $$\numlab{pimu}
It is easy to see that $\pi_\mu$ is a projection:
 $$
   \pi_\mu \circ \pi_\mu =
        \frac1{n^2} \sum_{j,l=0}^{n-1} \mu^{-j-l}
        \tilde\tau^{j+l} =
        \frac1n \sum_{l=0}^{n-1} \mu^{-l}\tilde\tau^l =
        \pi_\mu,
 $$
with $\im\pi_\mu \subset \Emu$:
 $$
   \tilde\tau\circ\pi_\mu = \frac1{n}
        \sum_{l=0}^{n-1} \mu^{-l} \tilde\tau^{l+1} =
        \mu \cdot \frac1{n}
        \sum_{l=0}^{n-1} \mu^{-l+1} \tilde\tau^{l+1} =
        \mu\cdot\pi_\mu.
 $$
We calculate that
 $$
   \sum_{\mu \in \zeta_n} \pi_\mu =
        \frac1n \sum_{\mu \in \zeta_n} \sum_{l=0}^{n-1}
        \mu^{-l} \tilde\tau^l
        = \frac1n \sum_{l=0}^{n-1}
        \Big(\sum_{\mu \in \zeta_n} \mu^{-l}\Big)
        \tilde\tau^l
         = \id_{\pi_*\E}.
 $$

It is clear from these properties of $\pi_\mu$ that
 $$
   \pi_*\E = \bigoplus_{\mu \in \zeta_n} \Emu. \tagn
 $$\numlab{disum1}
as sheaves of $\so_Y$-modules.

Let us now give an elementary argument for the fact that $\pi_*\E$
is locally free:

It is clear that $\E\big(\pi^{-1}(U)\big)$ for small enough $U$ is
a locally free $\so_X\big(\pi^{-1}(U)\big)$-module. Hence, to
establish that $\E\big(\pi^{-1}(U)\big)$ is a free
$\so_Y(U)$-module for small enough $U$, we just need to prove that
$\so_X\big(\pi^{-1}(U)\big)$ is a locally free $\so_Y(U)$-module
for small enough $U$. Away from the ramification points of $\pi$,
this is completely trivial. At a point $x\in \pi^{-1}(y)$ of a
special orbit $y\in P$, we consider a centered holomorphic
coordinate say $z$ such that $\tau^k z = \xi\cdot z$. It is easily
seen that $\pi_\mu(z^j)$, $j=0, \ldots, n'-1$, and $\mu\in\zeta_n$
such that $\xi^j=\mu^k$ provides a basis for
$\so_X\big(\pi^{-1}(U)\big)$ as an $\so_Y(U)$-module for small
enough $U$ around $y$. Here $\pi_\mu$ is defined just like in
(\numref{pimu}), except we use $\tau$ in place of $\tilde\tau$.

Using the fact that $\pi_*\E$ is locally free, equation
(\numref{disum1}) and the following exact sequence
 $$
   \xymatrix@C=0.6cm{
      \pi_*\E \ar[rr]^-{\oplus_{\mu'\neq \mu}\pi_{\mu'}} && \pi_*\E
      \ar[rr]^-{\pi_\mu} && \E_\mu \ar[r] & 0   },
 $$
we see that $\Emu$ is coherent and torsion free. Then $\Emu$ is
locally free, since $\dim_{\Bbb C}(Y) = 1$ (see Corollary (5.15)
in \cite{\refKo}). Let $E_\mu$ be the underlying holomorphic
bundle of $\Emu$. Then of course
 $$
   \pi_*E \cong \bigoplus_{\mu \in \zeta_n} E_\mu.
 $$
It is clear that $E_\mu$ has rank $2$ for all $\mu\in\zeta_n$.

We shall need the following technical lemma later, which can be
derived directly from the above.

\proclaim{Lemma \n\proclab{diag act}{Lemma}} Near each special
point $x\in X$ we can find a local frame $(s_1,s_2)$ for $E$ such
that $$
  \tilde\tau^k(s_i) = \theta_i \cdot s_i,
$$ where $\theta_i \in \zeta_{n'}$ are the eigenvalues of
$\tilde\tau^k$ acting on $E_x$. \endproclaim

\proof We observe, that if $\mu^k$ is different from $\theta_1$
and $\theta_2$ then any $s$ in the stalk $(\Emu)_x$
 must vanish at $x$. Now
choose $(s_\mu^{1},s_\mu^{2})$ local frame for $\Emu$ around $x$
for each $\mu \in \zeta_n$. Then $\{s_\mu^{i}\}$, $\mu \in
\zeta_n$ and $i=1,2$ is a local frame for $\pi_*E$. If
$\theta_1=\theta_2$ then this implies that
$\big(s_\mu^{1}(x),s_\mu^{2}(x)\big)$, where $\mu^k = \theta_1$,
is a basis of $E_x$, hence $(s_{\theta_1}^{1},s_{\theta_1}^{2})$
is the local frame we want in this case. If $\theta_1\neq\theta_2$
then there will be an $i_1$ and an $i_2$ such that
$\big(s_{\theta_1}^{i_1}(x),s_{\theta_2}^{i_2}(x)\big)$ is a basis
of $E_x$, hence $(s_{\theta_1}^{i_1},s_{\theta_2}^{i_2})$ is the
local frame we want in this case. \endproof

Our description of equivariant bundles and in turn the fixed point
set $|\M|$ is basically based on associating to an equivariant
bundle $(E,\tilde\tau)$ the bundle $E_1$ over Y together with the
numeric data $(d_1,d_2)$, but endowed with some extra structure,
namely a parabolic structure, since the bundle $E_1$ it self and
numeric data is not enough to determine the fixed point in
general. This parabolic structure is of course closely related to
the eigenspace decomposition of the fiber of $E$ over
$\pi^{-1}(P)$. To see how this goes, it is convenient to describe
how the bundle $E_1$ can be constructed using elementary
modifications, so let us here briefly recall the basics of
elementary modifications in this equivariant setting. --- In fact
we shall give a complete analysis of equivariant bundles, just
using these (inverse) elementary modification and only after that,
shall we as an aside return to the bundle $E_1$.

Let $y\in P$ be a special orbit, $(E,\tilde\tau)$ an equivariant
bundle of rank $2$ and $F_y = \oplus_{x\in \pi^{-1}(y)}F_x$, $F_x
\subset E_x$, be a $\tilde\tau$ invariant set of one-dimensional
flags over the orbit $y$. We then see that $F_x$ is an
eigensubspace for $\tilde\tau^k_x$. Let $\theta_1$ be the
eigenvalue of $\tilde\tau_x^k$ corresponding to the eigensubspace
$F_x$ and let $\theta_2$ be the other eigenvalue.

Let $S_x = E_x / F_x$ and let $S_y$ be the skyscraper sheaf on X
with support at $\pi^{-1}(y)$ and fiber $S_x$ at $x\in
\pi^{-1}(y)$.

\proclaim{Lemma \n\proclab{ex uniq el mod}{Lemma}} There is a
unique lift $(E',\tilde\tau')$ fitting into the following short
exact equivariant sequence
 $$
   \xymatrix@C=0.5cm{
     0 \ar[r] & E' \ar[rr]^-{\iota} && E \ar[rr]^-{\lambda} && S_y \ar[r] & 0
     },
 $$
such that $F_y = \oplus_{x\in\pi^{-1}(y)}\ker(\lambda_x) =
\oplus_{x\in\pi^{-1}(y)}\im(\iota_x)$. There is a set of flags
$\oplus_{x\in \pi^{-1}(y)} \ker(\iota_x)$ in $E'$ which consist of
eigensubspaces for $(\tilde\tau')_x^k$ whose corresponding
eigenvalue is $\xi^{-1}\theta_2$. The other eigenvalues of
$(\tilde\tau')_x^k$ is $\theta_1$. The determinants are related by
$\det E' \conge \det E \otimes [-\pi^{-1}(y)]$. \endproclaim

In this case we say that the equivariant bundle $(E',\tilde\tau')$
is obtained from the equivariant bundle $(E,\tilde\tau)$ by {\it
elementary modifications\/} in the direction $F_y$.

\proof Consider the sheaf kernel
 $$
   \E' = \ker\{\lambda : \E \longrightarrow S_y\}.
 $$
We observe, that $\E'$ is $\tilde\tau$-invariant. By using a local
frame adapted to the flag $F_x$ in $E_x$ one easily sees that
$\E'$ is locally free. The associated holomorphic vector bundle
$E'$ is uniquely determined, since $\E'$ as a subsheaf of $\E$ is
uniquely determined by $F_y$. Furthermore, since $\E'$ is
$\tilde\tau$-invariant, $\tilde\tau$ induces a lift $\tilde\tau'$
of $\tau$ to $E'$.

The determinant is clearly $\det E' \cong \det E \otimes (\det
S_y)^{-1}$, and $\det S_y$ can be calculated through the defining
sequence
 $$
   \xymatrix@C=0.5cm{
     0 \ar[r] & \so\big[-\pi^{-1}(y)\big] \ar[rr]^-{i} && \so_X \ar[rr]^-{p}
        && S_{y} \ar[r] & 0. } \tagn
 $$\numlab{trivsublb}
Again, using a frame of $E$ near $x$, which at $x$ gives an
eigenbasis for $\tilde\tau_x^{k}$ acting on $E_x$, one sees
immediately, that $ \ker(\iota_x)\subset E_x'$ is an eigenspace of
$(\tilde\tau')_x^k$ corresponding to the eigensubspace
$\xi^{-1}\theta_1$ and the other eigenvalue of $(\tilde\tau')_x^k$
is $\theta_2$. Using this it is easy to verify that the stated
isomorphism relation for the determinants also holds
equivariantly. \endproof

We observe, that the eigenvalue corresponding to $F_y$ is left
unchanged, where as the other eigenvalue is changed.

\remark{Remark \n\proclab{linesubs}{Remark}} We notice that if $L$
is an equivariant line sub-bundle of $E$ then
 $$
   L' = L \otimes \bigotimes_{\{x  \mid L_x \ne F_x\}}[-x]
 $$
is an equivariant line sub-bundle of $E'$. Conversely, if $L'$ is
a line sub bundle of $E'$ then $$
   L = L' \otimes \bigotimes_{\{x  \mid L'_x = \ker(\iota_x)\}}[x]
 $$
is a line sub-bundle of $E$ also equivariantly. This sets up a one
to one correspondence between line sub-bundles of $E$ and $E'$. By
considering the following commutative diagram of short exact
sequences for such a pair of corresponding equivariant line
sub-bundles $(L,L')$
 $$
   \xymatrix@R=0.3cm @!C=0.3cm{
     & 0 \ar[d] && 0 \ar[d]  && 0 \ar[d] \\
     0 \ar[r] & L' \ar[rr] \ar[dd] && L \ar[rr] \ar[dd]
        &&  S' \ar[dd] \ar[r] & 0 \\ \\
     0 \ar[r] & E' \ar[dd] \ar[rr]^{\iota} && E \ar[dd]
        \ar[rr]^{\lambda} && S_y \ar[dd] \ar[r] & 0 \\ \\
     0 \ar[r] & Q' \ar[d] \ar[rr] && Q \ar[d]
        \ar[rr] && S'' \ar[d] \ar[r] & 0 \\
       & 0 && 0  && 0 }
 $$
where $S' = \bigoplus_{ \{x \mid L_x \ne F_x\} }S_x$ and $S'' =
\bigoplus_{\{x \mid L_x = F_x\}}S_x$, we see there is a similar
equivariant correspondence between quotients and $$ Q' = Q \otimes
\bigotimes_{\{x \mid L_x = F_x\}}[-x]. $$

Going in the opposite direction we let $S'_x = T_x X \otimes F_x$
and $S'_y$ be the skyscraper sheaf on $X$ with support at
$\pi^{-1}(y)$ and fiber $S'_x$ at $x\in \pi^{-1}(y)$.

\proclaim{Lemma \n\proclab{ex uniq inv el mod}{Lemma}} There is a
unique lift $(E',\tilde\tau')$ fitting into the following short
exact equivariant sequence
 $$
   \xymatrix@C=0.5cm{
     0 \ar[r] & E \ar[rr]^-{\iota} && E'
        \ar[rr]^-{\lambda} && S'_y \ar[r] & 0
     },
 $$
such that $F_y = \oplus_{x\in\pi^{-1}(y)}\ker(\iota_x)$. There is
a set of flags $\oplus_{x\in \pi^{-1}(y)}\im(\iota_x)$ in $E'$
consisting of eigenspaces of $(\tilde\tau')_x^k$ whose
corresponding eigenvalue is $\theta_2$. The other eigenvalues of
$(\tilde\tau')_x^k$ is $\xi\theta_1$. \endproclaim

In this case we say that the equivariant bundle $(E',\tilde\tau')$
is obtained from the equivariant bundle $(E,\tilde\tau)$ by {\it
inverse elementary modifications\/} in the direction $F_y$.

\proof Let $\tilde\E$ be the sheaf of meromorphic sections of $E$
which are holomorphic everywhere except at $\pi^{-1}(y)$, where
the sections have a pole of order at most $1$. Consider now the
residue morphism $$ \res : \tilde\E \longrightarrow T_y X \otimes
E_y, $$ and compose it with the natural quotient map to $T_y X
\otimes S_y$ to obtain the composite morphism $\res_{F_y}$. Now
simply consider the sheaf kernel $$ \E' = \ker \{ \res_{F_y} :
\tilde\E \longrightarrow T_y X \otimes S_y\}. $$ By the very
construction of $\E'$, we see that it is $\tilde\tau$-invariant.
Again by considering a local frame adapted to $F_x$ it is easy to
see that $\E'$ is locally free. It clearly fits into the above
exact sequence and is uniquely determined by $F_y$ as a subsheaf
of $\tilde\E$. By the very construction we have that
$\ker(\iota_x) = F_x$ and by further assuming that the local frame
is an eigenbasis over $x$ one gets the statement about the
eigenvalues. \endproof

The statement we made about the determinants, line sub-bundles and
quotients in the case of elementary modification of course also
applies suitably adapted to inverse elementary modifications.

Let us now specify the iterations of (inverse) elementary
modifications, we shall use. Let $(E,\tilde\tau)$ an equivariant
rank $2$ bundle. We are interested in devising iterations of
elementary modification which at each of the special orbits $y\in
P$ changes one of the eigenvalues by $\xi^{-1}$ to the power, say
$m(y)\in \Z$ while keeping the other eigenvalue fixed. So, given
the {\it multiplicities\/} $m(y) \in \Z$ we proceed as follows
(under the convention, that if $m(y)$ is positive we will iterate
elementary modification at the orbit $y$ and if $m(y)$ is negative
we will iterate inverse elementary modification at $y$):

Suppose we are given a set of $\tilde\tau$-equivariant flags $F$
of $E|_{\pi^{-1}(P)}$. These flags and the signs of the
multiplicities $m$ then {\it determines\/} eigenvalues $\theta(y)$
at each orbit $y$ of non-zero multiplicity by the following rule:

At an orbits $y$, where $\tilde\tau_x^{k(y)}$ has two distinct
eigenvalues and $m(y)$ is positive, $F_y$ corresponds to the
eigenvalue $\theta(y)$ and where $m(y)$ is negative, $F_y$ is
complementary to the set of eigenspaces corresponding to the
eigenvalue $\theta(y)$.

We note that at the points $y\in P$, where the eigenvalues are
distinct, $\theta(y)$ and $\sign\big(m(y)\big)$ uniquely
determines $F_y = \oplus_{x\in \pi^{-1}(y)} F_x$. Applying either
elementary or inverse elementary modification, according to the
sign of $m(y)$, once to $(E,\tilde\tau)$ at $F_y$ for a $y\in P$,
results in a new equivariant bundle $(E',\tilde\tau')$. We can
define a set of flags $F'$ by the following assignments: \roster
    \item If $y'\in P -\{y\}$ then let $F'_{y'}= F_{y'}$.
    \item If $(\tilde\tau')_x^{k(y)}$ has two distinct eigenvalues for $x
        \in \pi^{-1}(y)$, then
        \vskip-\parskip
        \itemitem\llap{(a)\enspace}\ignorespaces if $m(y)$ is positive then let $F'_y$ be the
        set of flags given by the
        eigenspaces of $(\tilde\tau')_x^{k(y)}$, $x \in \pi^{-1}(y)$,
        corresponding to the eigenvalue $\theta(y)$.
        \vskip-\parskip
        \itemitem\llap{(b)\enspace}\ignorespaces if $m(y)$ is negative then let $F'_y$ be the
        set of flags given by the
        eigenspaces of $(\tilde\tau')_x^{k(y)}$, $x \in \pi^{-1}(y)$,
        corresponding to the eigenvalue different from $\theta(y)$.
    \item If $(\tilde\tau')_x^{k(y)}$ has only one eigenvalue for
    $x \in \pi^{-1}(y)$, then let $F'_y$ be the set of flags specified by
    the (inverse) elementary modification construction
    (see \procref{ex uniq el mod} and \noprocref{ex uniq inv el mod}).
\endroster

By the construction of $F'$, we see that it determines the same
eigenvalues $\theta(y)$, $y\in P$ as $F$ did. --- Note that under
\therosteritem2, $F'_y$ is set to be complementary to the set of
flags specified by the (inverse) elementary modification
construction.

It is the operation of taking $(E,\tilde\tau,F)$ to
$(E',\tilde\tau',F')$ through the (inverse) elementary
modification construction specified as above, we shall iterate
$|m(y)|$-times at each $y\in P$.

\definition{Definition \n\proclab{iterate i-elem mod}{Definition}}
Let $(E,\tilde\tau)$ be an equivariant rank $2$ bundle. Let $m: P
\to \Bbb Z$ be a multiplicity and $F$ be a set of
$\tilde\tau$-invariant flags for $E$ over $y\in P$ such that
$m(y)\neq 0$. We then define $\Gamma_{(m,F)}(E,\tilde\tau)$ to be
the equivariant rank $2$ holomorphic bundle obtained by iterating
elementary modification at $\{\, y\in P \mid m(y) \geq 0\, \}$
$m(y)$-times and inverse elementary modification at $\{\, y\in P
\mid m(y) < 0\, \}$ $|m(y)|$-times as described above.
\enddefinition

Note that there is a naturally induced set of flags in
$\Gamma_{(m,F)}(E,\tilde\tau)$ by its very construction.

\remark{Remark \n\proclab{og carlo}{Remark}} We observe, that the
equivariant bundle $\Gamma_{(m,F)}(E,\tilde \tau)$ together with
the induced set of flags, say $F'$, uniquely determines $(E,\tilde
\tau)$, since
 $$
   (E,\tilde\tau) = \Gamma_{(-m,F')}\big(\Gamma_{(m,F)}(E,\tilde \tau)\big).
 $$

\remark{Remark \n\proclab{og boev}{Remark}} Suppose
$(E,\tilde\tau)$ is an equivariant bundle with eigenvalues
$\theta_1 = \theta_1(y)$ and $\theta_2 = \theta_2(y)$ at the
special orbits $y\in P$. Let $m$ be a multiplicity and let $F$ be
 a set of flags which determines
$\theta_1$ as described above. Then the eigenvalues of the
equivariant bundle $(E',\tilde\tau') =
\Gamma_{(m,F)}(E,\tilde\tau)$ are $\theta_1(y)$ and
$\xi^{-m(y)}\theta_2(y)$. \endremark

This iteration will shortly be applied to equivariant bundles with
eigenvalues different from $1$, in order to associate other
equivariant bundles, which have all eigenvalues equal to $1$. The
significance of this is clear from the following lemma.

\proclaim{Lemma \n\proclab{equibndlpb}{Lemma}} An equivariant
bundle $(E,\tilde \tau)$ is equivariantly isomorphic to a pull
back bundle from $Y$ if and only if its numeric data vanishes.
\endproclaim

\proof It is clear that a pullback has this property. Conversely,
suppose $(E,\tilde\tau)$ has this property. We then claim that the
natural bundle map from $\pi^*(E_1)$ to $E$ is an equivariant
isomorphims. This is easily seen using the eigenframe for $E$
provided by \procref{diag act}. \endproof

We will further need the following lemma concerning equivariant
line bundles.

\proclaim{Lemma \n\proclab{invdivrep}{Lemma}} For any equivariant
line bundle $L$, there exists a $\tau$-invariant divisor $D$ such
that $[D]$ is equivariantly isomorphic to L. The difference
between any two such divisors is the pullback of a principal
divisor on $Y$. \endproclaim

\proof Let $\ttau$ be a lift of $\tau$ to $L$ and $\theta(y)$ be
the eigenvalues of $\ttau^{k(y)}$ acting on $L_x$, $x\in
\pi^{-1}(y)$, $y\in P$. Define $0\leq d(y) < n'(y)$, $y\in P$ by
$$
 \theta(y) = \xi(y)^{d(y)}
 $$
and set $D = \sum_{y\in P}d(y) \cdot \pi^{-1}(y)$. The induced
lift $\ttau'$ on the bundle $L' = L \otimes [-D]$ has all
eigenvalues equal to $1$. By arguing just as in the proof of the
\procref{equibndlpb}, it is easily seen that there exists a line
bundle $\bar L'$ on $Y$, such that $L' \conge \pi^*(\bar L')$. Let
$\bar D'$ be a divisor on $Y$ representing $\bar L'$. Then $D_L =
\pi^*(\bar D') + D$ is clearly a $\tau$-invariant divisor and $L
\conge [D_L]$.

Now note that for any $\tau$-invariant divisor $D'_L$ such that $L
\conge [D'_L]$, we must have $$ D'_L(y) = d(y) \mod n'(y), \text{
} \text{ for all } y \in P, $$ since the eigenvalues must agree.
But then $D'_L - D_L$ is clearly a pullback of a divisor from $Y$.
If $f\in {\Cal M}(X)$ such that $[D_L' - D_L] \conge (f)$, then
$f$ has to be $\tau$-invariant, since $[D_L' - D_L] \conge \so_X$.
Hence, there exits $g\in {\Cal M}(Y)$ such that $(\pi^*g)= (f) =
[D_L' - D_L]$. \endproof

\proclaim{Corollary \n\proclab{numcritpb}{Corollary}} For a
$\tau$-invariant divisor $D$, the {\it numeric data\/} $D(y) \mod
n'(y)$, $y\in P$, is an invariant of the isomorphism class of the
equivariant bundle $[D]$. \hfill $\square$ \endproclaim

If $L$ is an equivariant line bundle, we will write $0 \leq L(y) <
n'(y)$ for the numeric data $D(y) \mod n'(y)$, $y\in P$, where $D$
is any $\tau$-invariant divisor representing $L$ equivariantly.
Clearly, this numeric data of a line bundle is related to the
eigenvalues of the lift just like for rank $2$ bundles, but for
line bundles there is of course only one eigenvalue per point in
$P$.

\proclaim{Corollary \n\proclab{numcritpb}{Corollary}} An
equivariant line bundle $L$ is equivariantly isomorphic to a pull
back bundle from $Y$ if and only if its numeric data $L(y)$
vanishes at every $y\in P$.\hfill $\square$ \endproclaim

Let us now return to the rank $2$ situation, so let $(E, \ttau)$
be an equivariant rank $2$ bundle with numeric data $(d_1,d_2)$
and corresponding eigenvalues $(\theta_1,\theta_2)$. Let $m(y) =
d_2(y) - d_1(y)$ and define a $\tau$-invariant divisor $D_2$ on
$X$ by $D_2 = \sum_{y\in P} d_2(y) \cdot \pi^{-1}(y)$. Consider
the holomorphic bundle $E' = E \otimes [-D_2]$ with the induced
lift $\tilde\tau'$. The eigenvalues of $(\tilde\tau')^k$ are $1$
and $\theta(y) = \theta_1(y)\theta_2(y)^{-1}$.

For $x\in \pi^{-1}(P)$, such that $m\big(\pi(x)\big)\neq 0$, we
let $F'_x \subset E'_x$ be the eigenspace of $(\tilde\tau')_x^k$
corresponding to the eigenvalue $\theta(y)$ and let $F'_y =
\oplus_{x\in \pi^{-1}(y)}F'_x$.

Observe, that the equivariant bundle
$\Gamma_{(-m,F')}(E',\tilde\tau')$ induces the identity on fibers
over every $x \in \pi^{-1}(P)$, hence there is a unique
quasi-parabolic bundle $(\bar E, \bar F)$ on $Y$ such that
$\big(\pi^*\bar E, \pi^*\bar F\big)$ with the naturally induced
pullback lift is equivariant isomorphic to
$\Gamma_{(-m,F')}(E',\tilde\tau')$, and such that the flags
$\pi^*\bar F$ gets identified with the induced flags of
$\Gamma_{(-m,F')}(E',\tilde\tau')$.

We give $\big(\bar E, \bar F\big)$ the structure of a parabolic
bundle by letting the weights $w(y)$ associated to the parabolic
points $y\in P$, $m(y)\neq 0$, be given by
 $$
   w(y) = \frac{m(y)}{n'(y)}.\tagn
 $$\numlab{par weight}
Denote the resulting parabolic bundle by $(\bar E, \bar F, w)$. In
the notation of \cite{\refSe}, $w(y) = a_{y2} -a_{y1}$.

We observe, that the determinants of $E$ and of $\bar E$ are
equivariantly related by
 $$
   \pi^*\det(\bar E) \cong \det(E) \otimes [D-2D_2],
 $$
where $D = \sum_{y\in P}m(y) \cdot \pi^{-1}(y)$.

Let $\Delta \in \Dc$ be the equivariant determinant of
$(E,\ttau)$. Then from \procref{numcritpb} we get the following
conditions on the numeric invariants of $(E,\ttau)$: $$
 \Delta(y) = d_1(y) + d_2(y) \mod n'(y),\tagn
$$\numlab{resnd} for all $y\in P$. Define the a line bundle
$\bar\Delta$ on $Y$ by requiring $$ \pi^*\bar\Delta \conge \Delta
\otimes \big[-\sum_{y\in P} \big(d_1(y) + d_2(y)\big) \cdot
\pi^{-1}(y)\big].\tagn $$\numlab{bardelta} Clearly, there is a
unique such line bundle on $Y$ up to isomorphism and $$
  \det(\bar E)\cong \bar\Delta.
$$

\proclaim{Lemma \n\proclab{parbndluniquedetfp}{Lemma}} Given $\Cal
D$, the isomorphism class of the parabolic bundle $(\bar E, \bar
F,w)$ together with $d_2(y)$, $y\in P$ uniquely determines the
equivariant bundle $(E,\tilde\tau)$. In the case $n$ is odd, the
parabolic bundle $(\bar E, \bar F,w)$ alone determines the
equivariant bundle $(E,\tilde\tau)$. In case $n$ is even, the
parabolic bundle $(\bar E, \bar F,w)$ uniquely determines
$\Delta(y) - 2d_2(y) \mod n'(y)$ for all $y\in P$. \endproclaim

Whether $(\bar E, \bar F,w)$ actually also determines
$(E,\tilde\tau)$ or not in case $n$ is even, comes down to some
delicate questions about existence of meromorphic functions with
divisors of certain special kind with support contained in
$\pi^{-1}(P)$. We leave it to the reader to check the general
fact, that if $n'(y)$ is even, for one of the $y\in P$, then
$\Delta \in \Dc$ is also determined by $(\bar E, \bar F,w)$.

\proof That $(\bar E,\bar F,w)$ together with $d_1(y)$, $y\in P$
determines $(E,\ttau)$ is immediate from \procref{og carlo}. The
determinant relation, which holds for some $\Delta\in \Dc$ states
that,
 $$
  \pi^* \det \bar E \cong \Delta \otimes [D - 2 D_2].
 $$
Hence,
 $$
  0 = \Delta(y) + n'(y)w(y) - 2 d_2(y) \mod n'(y).
 $$
In case $n$ is odd, $\Delta(y)$ is uniquely determined from the
outset and since $n'(y)$ must be odd in this case, we see $d_2$ is
uniquely determined from this.

In case $n$ is even, the wanted conclusion follows immediately
from this equation. \endproof

Let us now reverse the process and provide a construction of all
equivariant bundles
 from a certain set of admissible parabolic bundles. Clearly the
 numeric data need to satisfy equation (\numref{resnd}), so we
 make the following definition.

\definition{Definition \n\proclab{labelset}{Definition}}
For each $\Delta \in \Dc$, define the subset $$
  \LD \subset \prod_{y\in P}T_{n'(y)}
$$ by $(d_1,d_2)\in \LD$ if and only if $$ \Delta(y) = d_1(y) +
d_2(y) \mod n'(y), $$ for all $y\in P$. In case $P=\emptyset$, we
let $\LD = \{(0,0)\}$. For each pair $(d_1,d_2)\in\LD$ we let the
line bundle $\bar \Delta$ be given by formula (\numref{bardelta})
and weights $w : P \ra [0,1)$ given by $$ w(y) = \frac{d_2(y) -
d_1(y)}{n'(y)}, $$ for all $y\in P$. \enddefinition

We note that $\bar \Delta$ and $w$ are uniquely determined by
$\Delta$ and the pair $(d_1,d_2) \in\LD$.

Let $\tilde P (\bar\Delta,w)$ be the set of isomorphism classes of
rank $2$ parabolic bundles on $(Y,P)$ with the parabolic weight
$w(y)$ at $y\in P$ and determinant isomorphic to $\bar\Delta$ and
let $P(\bar\Delta,w)$ be the moduli space of semi-stable parabolic
bundles on $(Y,P)$ with the parabolic weight $w(y)$ at $y\in P$
and determinant isomorphic to $\bar\Delta$. We note here that no
quasi-parabolic structure is present at a point $y$, if $w(y)=0$.

\definition{Definition \n\proclab{admisparabb}{Definition}}
The {\it set of admissible parabolic bundles\/} on $Y$ is defined
to be $$ \tilde P_a = \coprod_{\Delta\in\Dc}\coprod_{(d_1,d_2)\in
\LD} \tilde P(\bar\Delta,w), $$ and the {\it moduli space of
admissible parabolic bundles\/} on $Y$ is defined to be $$ P_a =
\coprod_{\Delta\in\Dc}\coprod_{(d_1,d_2)\in \LD} P(\bar\Delta,w).
$$ \enddefinition

Now let $(d_1,d_2)\in \LD$, for $\Delta\in \Dc$ and suppose $(\bar
E,\bar F,w)$ is a parabolic bundle representing an element in
$\tilde P(\bar\Delta,w)$. Then we define
 $$
  (E,\tilde\tau) = \Gamma_{(wn',\pi^*\bar F)}(\pi^*\bar
  E)\otimes[D_2],\tagn
 $$\numlab{defF}
where $D_2 = \sum_{y\in P}d_2(y) \cdot \pi^{-1}(y)$ as before and
$\pi^*(\bar E)$ is given the pullback lift.

\proclaim{Theorem \n\proclab{bijection}{Theorem}} Let $\tilde\F :
\tilde P_a \ra \Li$ be the map given by {\rm(\numref{defF})}
above. Then $\tilde \F$ is a bijection and the inverse map is
given by the construction described just after
\procref{numcritpb}. \endproclaim

\proof Immediate from \procref{og carlo} \endproof

We will denote the inverse map of $\tilde \F$ by $\tP : \Li \ra
\Pa$.

Let us end this section by establishing the relation between the
underlying bundle of $\tP(E,\ttau)$ and $E_1$, for any equivariant
bundle $(E,\ttau)$.

\proclaim{Lemma \n\proclab{Eone}{Lemma}} Suppose $(E,\ttau)$ is an
equivariant bundle and $(\bar E,\bar F, w) = \tP(E,\ttau)$. Then
$\bar E \cong E_1$. \endproclaim

\proof Let $(d_1,d_2)$ be the numeric data of $(E,\ttau)$. Then it
is rather easy to see that
 $$
  \pi^*E_1 \conge \Gamma_{(m_2,F_1)}\circ \Gamma_{(m_1,F_2)}(E,\ttau)
          = \Gamma_{(m,F_1)}(E,\tilde\tau) \otimes [-D_1],
 $$
where $m_\nu(y) = d_\nu(y)$, $y\in P$, $F_\nu$ is the set of
eigenspaces corresponding to the eigenvalue $\theta_\nu$ over $P$,
and $D_1 = \sum_{y \in P} d_1(y) \cdot \pi^{-1}(y)$. From the
construction of $\bar E$, we have that
 $$
  \align
   \pi^*\bar E &= \Gamma_{(-m,F')}\big((E,\tilde\tau) \otimes [-D_2]\big)
                     = \Gamma_{(-m,F')} \circ \Gamma_{(m,G_2)}
                    \circ \Gamma_{(m,G_1)}\big((E,\tilde\tau) \otimes
                    [-D_1]\big) \\
        &= \Gamma_{(m,G_1)}\big((E,\tilde\tau) \otimes [-D_1]\big) =
                    \Gamma_{(m,F_1)}(E,\tilde\tau) \otimes [-D_1],
  \endalign
 $$
where $F'$ corresponds to the eigenvalue $\theta_1\theta_2^{-1}$,
$G_1$ to $1$ and $G_2$ is induced in
$\Gamma_{(m,G_1)}\big((E,\tilde\tau) \otimes [-D_1]\big)$ in the
usual fashion. This implies the lemma. \endproof

\head The fixed Points in the Moduli Space.\endhead \seclab{fixed
ge}

Recall that by \procref{repfixeb} any fixed point in $\M$ can be
represented by an equivariant rank $2$ bundle $(E,\tilde\tau)\in
\Li$ and that \procref{bijection} gives a complete classification
of equivariant bundles in terms of admissible parabolic bundles.

Once we have established that $E = \Pi_e\tF(\bar E,\bar F,w)$ is
semi-stable provided $(\bar E,\bar F,w)$ is parabolically
semi-stable and that S-equivalent semi-stable parabolic bundles go
to S-equivalent bundles under $\Pi_e\tF$, we get an induced
 map from $P_a$ to $|\M|$. Let us first address the
semi-stability issue.

\proclaim{Proposition \n\proclab{ssiffpss}{Proposition}} The
bundle $E$ is semi-stable if and only if $(\bar E,\bar F,w)$ is
parabolically semi-stable. \endproclaim

Based on this proposition, we see that if $\tPa'$ is defined to be
the subset of $\tPa$ consisting of isomorphism class of those
parabolic bundles, which are semi-stable, then we get an induced
bijection $\tF : \tPa' \ra \Li'$.

We need the following lemma to prove the Proposition.

\proclaim{Lemma \n\proclab{stab eq}{Lemma}} Suppose $L$ is a
$\tilde\tau$-invariant subbundle of $E$. Then there is an induced
parabolic subbundle $\bar L$ of $\bar E$ such that
 $$
   \pa\mu(\bar E) - \pa\mu(\bar L)
        = \tfrac1n \big(\mu(E)- \mu(L)\big).
 $$
Conversely, a parabolic subbundle $\bar L$ of $\bar E$ induces an
invariant subbundle $L$ of $E$, and the same equality holds.
\endproclaim

\proof From equation (\numref{defF}) relating $E$ and $\bar E$ it
is clear that there is a one to one correspondence between line
sub-bundles of $\bar E$ and $\tilde\tau$-invariant sub-bundles of
$E$. Let $(\bar L,L)$ be such a pair of corresponding line
sub-bundles.

Recall that when $\bar L$ is given the induced parabolic structure
from $(\bar E,\bar F,w)$, then we have that
 $$
  \pa\mu(\bar E) - \pa\mu(\bar L) = \tfrac12\deg \bar E - \deg \bar L
  + \tfrac12 \sum_{\bar L_y \neq \bar F_y} w(y) - \tfrac12
  \sum_{\bar L_y = \bar F_y} w(y).
 $$
Now
 $$
   \deg E = n \deg\bar E - \deg D + 2 \deg D_2,
 $$
and
 $$
   \deg L = n \deg\bar L - \sum_{\bar L_y \ne \bar F_y}kn'(y) w(y) + \deg D_2.
 $$
Computing $\frac12 \deg E - \deg L$ from the equations for $\deg
E$ and $\deg L$ and comparing with the above expression for
$\pa\mu(\bar E) - \pa \mu(\bar L)$, we get the formula we wanted.
\endproof

\demo{Proof of \procref{ssiffpss}} Suppose $(\bar E, \bar F, w)$
is parabolic semi-stable and let $W$ be the unique maximal
semi-stable sub-bundle of $E$. As $\tilde\tau(W)$ is a bundle with
the same properties, $\tilde\tau(W) = W$ by uniqueness. Assume
that $W \ne E$ so that $E$ is not semi-stable. Then $W$ is an
invariant line bundle in $E$ and employing the previous lemma we
get: $\mu(W) \le \mu(E)$ by semi-stability of $(\bar E, \bar F,
w)$, thus contradicting the assumption on $E$. Hence, we have
proved that $E$ is semi-stable whenever $(\bar E, \bar F, w)$ is
parabolic semi-stable.

If on the other hand $E$ is stable (resp.\ semi-stable), then it
follows immediately from \procref{stab eq} that $(\bar E, \bar F,
w)$ is parabolically stable (resp.\ semi-stable). \endproof

Let us now check that parabolic S-equivalence classes are taken to
S-equivalence class of semi-stable bundles:
$\Gr\big(\Gamma_{(wn',\pi^*\bar F)}(\bar E, \bar F,w)\big) =
\Gamma_{(wn',\pi^*\bar F)}\big(\Gr(\bar E, \bar F,w)\big)$. What
remains to be shown is that this is the case for semi-stable, but
non-stable bundles. So suppose that
 $$
   \xymatrix@C=0.5cm{
     0 \ar[r] & (\bar L_1,w_1) \ar[rr] && (\bar E, \bar F, w)
            \ar[rr] && (\bar L_2,w_2) \ar[r] & 0 },
 $$
is an exact sequence of semi-stable parabolic bundles such that
$\pa\mu (\bar L_\nu,w_\nu) = \pa\mu (\bar E, \bar F, w)$. Then
 $$
   \Gamma_{(wn',\pi^*\bar F)}\big(\pi^*(\bar L_1,w_1) \oplus \pi^*(\bar
      L_2,w_2)\big)  = (L_1,\tilde\tau_1)\oplus (L_2,\tilde\tau_2),
 $$
where according to \procref{linesubs} $(L_\nu, \tilde\tau_\nu) =
\pi^*\bar L_\nu \otimes [-\sum_{y\in P}w_\nu(y)n'(y)\cdot
\pi^{-1}(y)]$.

On the other hand, $\bar L_1$ gives a $\tilde\tau$-invariant
destabilising line subbundle $E$ which is clearly equivariantly
isomorphic to $(L_1,\tilde\tau_1)$. Put $(L'_2, \tilde\tau'_2) =
(E,\tilde\tau)/L_1$ so that $\Gr(E,\tilde\tau) =
(L_1,\tilde\tau_1) \oplus (L'_2, \tilde\tau'_2)$. Then it follows
from \procref{linesubs} that $(L'_2, \tilde\tau'_2) \conge (L_2,
\tilde\tau_2)$ and we get that $\Gr(\bar E,\bar F, w) = (\bar
L_1,w_1)\oplus (\bar L_2,w_2)$.

We have now established that the following map is well-defined:

\definition{Definition \n\proclab{DEFF}{Definition}}
We define a map $$\F : P_a \lra |\M|$$ by requiring commutativity
of $$
   \CD
     \tilde P_a' @>\tF>> \Li' \\
     @VVV @V VV  \\
     P_a @>\F>> |\M|.
   \endCD
 $$
\enddefinition

From the above discussion, it is clear that $\F$ is surjective and
we shall now discuss the fibers of $\F$. First let
 $$
   P_a = P_a^{ss}\sqcup P_a^{s},
 $$
where $P_a^{s}$ is the set of point in $P_a$ represented by stable
parabolic bundles and $P_a^{ss}$ is the complement of these in
$P_a$.

Let $T$ be the sub-group generated by the involutions in $\ker
\pi^* : \Pic(Y) \ra \Pic(X)$. If $r$ is even (which is not the
case if n is odd), then $T$ is the unique sub-group of this kernel
of order two. Otherwise, it is the trivial group $\{1\}$. In case
$r$ is even, $T$ will act on $P_a^{s}$ and we decompose
 $$
   P_a^{s} = P_a^{s,i} \sqcup P_a^{s,g},
 $$
where $P_a^{s,i} = (P_a^{s})^T$ is the fixed point set under the
$T$-action. If $r$ is odd, we let $P_a^{s,g} = P_a^{s}$ and
$P_a^{s,i} = \emptyset$.

Likewise, let $\M^s\subseteq M$ be the subset represented by
stable bundles and $\M^{ss} = \M - \M^s$.

Our main set-theoretic theorem states:

\proclaim{Theorem \n\proclab{mainsetT}{Theorem}} The map $\F : P_a
\ra |\M|$ is onto and has finite fibers.

In case $n$ is {\rm odd}
 $$
   \F : P_a^s \lra |\M^s|
 $$
is a bijection and
 $$
   \F : P_a^{ss} \lra |\M^{ss}|
 $$
is invariant under the finite\footnote{We call an equivalence
relation {\it finite\/}, if it has finite equivalence classes.}
equivalence relation $\sim_o$ described below under point 1.\ and
it induces a bijection between $P_a^{ss}/\sim_o$ and $|\M^{ss}|$.

In case $n$ is {\rm even}, $\F$ is invariant under the
$\zeta_2$-action described below under point 2.\ and it descends
to
 $$
   \F : \bar P_a = P_a/\zeta_2 \lra |\M|.
 $$
This $\zeta_2$ action preserves stability and coincides with the
$T$-action whenever $r$ is even and we get an induced
decomposition $\bar P_a^s = \bar P_a^{s,i} \coprod \bar
P_a^{s,g}$. The restriction
 $$
   \F : \bar P_a^{s,g} \lra |\M^s|
 $$
is a bijection and
 $$
   \F : \bar P_a^{s,i} \sqcup \bar P_a^{ss} \lra |\M^{ss}|
 $$
is invariant under the finite equivalence relation $\sim_e$ on
$\bar P_a^{ss}$ described under point 1.\ below and induces a
bijection between $\bar P_a^{s,i} \sqcup (\bar P_a^{ss}/\sim_e)$
and $|\M^{ss}|$. \endproclaim

\subhead 1. The equivalence relations $\sim_o$ and
$\sim_e$\endsubhead Using the $\zeta_n$-action on the set of
equivariant line bundles we introduce the following equivalence
relation.

\definition{Definition \n\proclab{DEFeqo}{Definition}}
We define the equivalence relation $\sim_o$ on $\Pa^{ss}$ by
declaring that $$ \Gr(\bar E,\bar F, w) \sim_o \Gr(\bar E',\bar
F', w') $$ if and only if there exists $\mu\in\zeta_n$ and
equivariant line bundles $(L_\nu,\ttau_\nu)$, $\nu=1,2$, such that
$$
  \gather
   \tF\big(\Gr(\bar E,\bar F, w)\big) \cong (L_1,\ttau_1) \oplus (L_2,\ttau_2), \\
   \tF\big(\Gr(\bar E',\bar F', w')\big) \cong
        (L_1,\mu\cdot\ttau_1) \oplus (L_2,\mu^{-1}\cdot\ttau_2).
  \endgather
 $$
\enddefinition

Let us now examine how the numeric data $(d_1,d_2)$ of
$\tF\big(\Gr(\bar E,\bar F, w)\big)$ is related to the numeric
data $(d_1',d_2')$ of $\tF\big(\Gr(\bar E',\bar F', w')\big)$. Let
$f_\mu$ be a meromorphic function on $X$, such that $f_\mu \circ
\tau^{-1} = \mu\cdot f_\mu$ and let $d_\mu(y) = (f_\mu)(y)$, $y\in
P$. Further let $\Gr(\bar E,\bar F, w) = (\bar L_1,w_1) \oplus
(\bar L_2,w_2)$ and $\Gr(\bar E',\bar F', w') = (\bar L_1',w_1')
\oplus (\bar L_2',w_2')$, where
 $$
   w_\nu(y) = \cases 0, &\quad\text{for } \bar L_{\nu,y} \ne \bar F_y, \\
            w(y) &\quad\text{for } \bar L_{\nu,y} = \bar F_y \endcases
 $$
and $w'_\nu$ similarly. Set $0\leq d^{\nu}(y) < n'(y)$, $y\in P$,
to be $$
  \gather
   d^{1}(y) = d_2(y) - n'(y)w_1(y) + d_\mu(y) \mod n'(y), \\
   d^{2}(y) = d_2(y) - n'(y)w_2(y) - d_\mu(y) \mod n'(y).
  \endgather
 $$
Then if $P_\times = \{\, y\in P \mid d^{1}(y) > d^{2}(y)\, \}$, we
have that $$ \big(d_1', d_2'\big)(y) = \cases
        \big(d^{1}(y), d^{2}(y)\big), &\quad\text{for } y\in P- P_\times, \\
        \big(d^{2}(y), d^{1}(y)\big), &\quad\text{for } y\in
        P_\times.
   \endcases\tagn
$$\numlab{Ldeqo} Of course $w'(y) =
\big(d_2'(y)-d_1'(y)\big)\big/n'(y)$, $w'_\nu(y) =
d^{\nu}(y)/n'(y)$ and $\bar F'_y = \bar L'_{\nu,y}$ if and only if
$\bar F_y = \bar L_{\nu,y}$ for $y \in P - P_\times$ and for $y
\in P_\times$ the flag is switched: $\bar F'_y = \bar
L'_{\nu+1,y}$ if and only if $\bar F_y = \bar L_{\nu,y}$.

\proclaim{Lemma \n\proclab{simo}{Lemma}} The numeric data
$(d_1,d_2)$ together with the subset $P_1 = \{\, y\in P \mid
w_1(y) \neq 0\,\}\subseteq P$ and $\mu$ determines a line bundle
$L$ such that $$
  \bar L_1' \cong \bar L_1 \otimes L, \quad\text{and}\quad
   \bar L_2' \cong \bar L_2 \otimes L^{-1}.
$$ \endproclaim

\proof Let $\bar D_\nu$ be divisors on $Y$ such that $\bar L_\nu =
[\bar D_\nu]$. Then let
 $$
  D^1 = \pi^*\bar D_1 + D_2 + (f_\mu) - \sum_{y\in P} n'(y) w_1(y)
  \cdot \pi^{-1}(y).
 $$
By the formula, we see that $[D^1] \cong (L_1,\mu\cdot\ttau_1)$ as
an equivariant bundle. We now observe that there is a divisor
$\bar D_1'$ on $Y$ uniquely determined by
 $$
  \pi^*\bar D_1' = D^1 - D'_2  + \sum_{y\in P} n'(y) w'_1(y)
  \cdot \pi^{-1}(y),
 $$
where $D'_2 = \sum_{y \in P}d'_2(y) \cdot \pi^{-1}(y)$. From this
we see that $\pi^*\bar D_1' - \pi^*\bar D_1 = \pi^*\bar D$ and
$\bar D$ is determined completely by $(d_1,d_2)$ and $P_1$. By
calculating the determinants, one can check that $L = [\bar D]$
has the required properties. \endproof

In case $n$ is {\it even\/}, we define an equivalence relation
$\sim_e$ on $\Pa^{ss}$ generated by $\sim_o$ just described and
then the following equivalence relation $\sim$: Let $\mu\in
\zeta_n$ be such that $\mu\cdot \Delta \cong \Delta'$ as
equivariant bundles, where $\Dc = \{\Delta,\Delta'\}$.

\definition{Definition \n\proclab{DEFeqe}{Definition}} The $\Gr(\bar E,\bar F, w)\in
\Pa^{ss}(\Delta)$ is $\sim$-equivalent to $\Gr(\bar E',\bar F',
w') \in \Pa^{ss}(\Delta')$ (written $\Gr(\bar E',\bar F', w') \sim
\Gr(\bar E,\bar F, w)$) if and only if there exists equivariant
line bundles $(L_\nu,\ttau_\nu)$, $\nu=1,2$, such that $$
  \gather
   \tF\big(\Gr(\bar E,\bar F, w)\big) \cong (L_1,\ttau_1) \oplus (L_2,\ttau_2), \\
   \tF\big(\Gr(\bar E',\bar F', w')\big) \cong (L_1,\mu\cdot\ttau_1) \oplus (L_2,\ttau_2).
  \endgather
 $$
\enddefinition Again, we let $\Gr(\bar E,\bar F, w) = (\bar
L_1,w_1) \oplus (\bar L_2,w_2)$ and $\Gr(\bar E',\bar F', w') =
(\bar L_1',w_1') \oplus (\bar L_2',w_2')$. Set $0\leq d^{\nu}(y) <
n'(y)$, $y\in P$, to be $$
  \gather
   d^{1}(y) = d_2(y) - n'(y)w_1(y) + d_\mu(y) \mod n'(y), \\
   d^{2}(y) = d_2(y) - n'(y)w_2(y)  \mod n'(y).
  \endgather
 $$
Then $d_\nu'$ and $w_\nu'$ are given by (\numref{Ldeqo}) and the
lines following it.

\proclaim{Lemma \n\proclab{sime}{Lemma}} The indexing of the line
bundles $\bar L_\nu$ and $\bar L_\nu'$, $\nu=1,2$, can be chosen
such that $$
   \bar L_1' \cong \bar L_1 \otimes \bar \Delta_c,
   \quad\text{and}\quad
   \bar L_2' \cong \bar L_2,
 $$
where $\bar \Delta_c = \bar \Delta' \bar \Delta^{-1}$.
\endproclaim

\proof Obvious. \endproof

\subhead 2. The action of $\zeta_2$ on $P_a$\endsubhead In the
case where $n$ is even, there is an natural $\zeta_2$-action on
$\Li$, simply gotten by letting $-1\in \zeta_2$ map an equivariant
bundle $(E,\tilde\tau)$ to the equivariant bundle $(E,-\ttau)$,
which of course maps $\Li'$ to it self. Using $\tF$, we get an
induced $\zeta_2$ action on $\tPa'$, which by the properties of
$\tF$ induces an action on $\Pa$. Let us now explicitly describe
this action on $\Pa$.

If $(d_1,d_2)$ is the numeric data of $\tF(\bar E,\bar
F,w)=(E,\tilde\tau)$ and $(d'_1,d'_2)$ that of $\tF(\bar E',\bar
F',w')=(E,-\ttau)$. We observe, that the numeric data are related
by
  $$
   \big(d_1', d_2'\big)(y) = \cases
        \big(d_1(y), d_2(y)\big), &\quad\text{for } y \in P-P_o, \\
        \big(d_1(y) + \tfrac {n'(y)}2, d_2(y) + \tfrac {n'(y)}2\big), &\quad\text{for } y \in P_+, \\
        \big(d_1(y) - \tfrac {n'(y)}2, d_2(y) - \tfrac {n'(y)}2\big), &\quad\text{for } y \in P_-, \\
        \big(d_2(y) - \tfrac {n'(y)}2, d_1(y) + \tfrac {n'(y)}2\big), &\quad\text{for } y \in
        P_\times,
   \endcases\tagn
 $$\numlab{LDvirk2}
where
 $$
  \gather
   P_o = \{\, y \in P \mid k(y) \text{ odd} \, \}, \qquad
   P_+ = \{\, y \in P_o \mid d_2(y) < \tfrac{n'(y)}2 \, \}, \\
   P_- = \{\, y \in P_o \mid d_1(y) \ge \tfrac{n'(y)}2 \, \},
   \qquad
   P_\times = \{\, y \in P_o \mid d_1(y) < \tfrac{n'(y)}2 \le d_2(y) \,
   \}.
  \endgather
 $$
Let $m_\times(y) =1$ for $y\in P_\times$ and zero otherwise.
Clearly, both $m_\times$ and $w'$ are determined just by
$(d_1,d_2)$.

\proclaim{Lemma \n\proclab{act2}{Lemma}} The numeric data
$(d_1,d_2)$ determines a line bundle $\bar L$ on $Y$, such that $$
(\bar E',\bar F') = \Gamma_{(-m_\times,\bar F)}(\bar E, \bar F)
\otimes \bar L.\tagn $$\numlab{act2f} \endproclaim

Observe, that if $r$ is even then all $k(y)$ are even and
therefore $m_\times(y)=0$ for all $y\in P$, so $(\bar E',\bar F')
= (\bar E, \bar F) \otimes \bar L$. But then $\bar L\in \ker
\pi^*$ and in fact $\pi^*(\bar L) \conge [(f_{-1})]$, which means
that $T=\langle \bar L \rangle$.

\proof Let $D_- = (f_{-1})$. Then we have an isomorphism of
equivariant bundles
 $$
  \Gamma_{(wn',\pi^*\bar F)}(\pi^*\bar E) \otimes [D_2] \otimes
      [D_-] \conge (E,-\ttau) \conge \Gamma_{(w'n',\pi^*\bar
      F')}(\pi^*\bar E') \otimes [D_2'].
 $$
Let $F'$ be the flag induced in $E$ from this isomorphism. Then
 $$
  \pi^*\bar E' \conge \Gamma_{(-w'n',F')}\circ\Gamma_{wn',\pi^*\bar F}(\pi^*\bar
          E) \otimes [D_- - D_2' + D_2].
 $$
Now
 $$
  w'(y) = \cases
          w(y), &\quad\text{for } y \in P-P_\times \\
          1-w(y), &\quad\text{for } y \in P_\times,
     \endcases
 $$
so at $y\in P -P_\times$ no modification is done to $\pi^*\bar E$,
but at each $y\in P_\times$ we do $-n'(y)$ elementary
modifications in the direction $\pi^* \bar F_y$, followed by
$n'(y)w(y)$ elementary modifications in two complementary
directions, which just corresponds to tensoring with $[-n'(y)w(y)
\cdot \pi^{-1}(y)]$. Hence, we see that
 $$
  \pi^*\bar E' \conge \Gamma_{(-n'm_\times,\pi^*\bar F)}(\pi^*\bar E)
  \otimes [-D_\times + D_- - D_2' + D_2],
 $$
where $D_\times = \sum_{y \in P_\times}w(y) n'(y) \cdot
\pi^{-1}(y)$. But now observe that $(-D_\times + D_- - D_2' +
D_2)(y)$ is zero $\mod n'(y)$ for all $y\in P$, hence $(-D_\times
+ D_- - D_2' + D_2)$ is a pullback of a divisor, say $\bar D$,
from $Y$ and we get that $$ \bar E' \cong \Gamma_{(-m_\times,\bar
F)}(\bar E) \otimes [\bar D]. $$ \endproof

Note that this $\zeta_2$ action on $P_a^{ss}$ is compatible with
the equivalence relation $\sim_e$, thus we get an induced
equivalence relation of $\bar P_a^{ss}$, which we also denote
$\sim_e$.

In order to prove \procref{mainsetT}, we need to understand which
stable parabolic bundles in $\Pa^s$ are mapped to non-stable (but
of course semi-stable) bundles under $\F$:

\proclaim{Lemma \n\proclab{char of deg}{Lemma}} Let $(\bar E,\bar
F,w)\in \Pa^{s}$. Then $\F(\bar E, \bar F,w) \in |\M^{ss}|$ if and
only if $(\bar E,\bar F,w)\in \Pa^{s,i}$, i.e. $r$ is even and
$(\bar E,\bar F,w) \cong (\bar E,\bar F,w) \otimes \bar L$, where
$T = \langle \bar L\rangle$. \endproclaim

\proof Assume that $\tF(\bar E,\bar F,w) = (E,\ttau)$ and
$\Gr(E)=L_1\oplus L_2\in |\M^{ss}|$, where $L_1\subset E$. Since
$\tau^*\Gr(E) \cong \Gr(E)$, we have either that $\tau^*L_1 \cong
L_1$ and $\tau^* L_2 \cong L_2$ or $\tau^* L_1 \cong L_2 \ncong
L_1$.

First the invariant case, $\tau^*L_\nu \cong L_\nu$. Consider the
sub-case, where $\ttau(L_1) = L_1$. Then this $\ttau$-invariant
subbundle $L_1$ induces a parabolic subbundle $(\bar L_1,w_1)
\subset (\bar E,\bar F,w)$ with $\pa\mu(\bar L_1,w_1) =
\pa\mu(\bar E,\bar F,w)$ simply by following $L_1$ through the
construction of $(\bar E,\bar F,w)$ from $(E,\ttau)$. But this
contradicts the stability of $(\bar E,\bar F,w)$. Now, in the
sub-case where $\ttau(L_1)\neq L_1$ then $L_1\cap \ttau(L_1) =
\{0\}$ and $\ttau (L_1) \cong L_2$, so $E \cong L_1 \oplus
\ttau(L_1)$ and $L_1 \cong L_2$, i.e. $E\cong L_1 \oplus L_1$. It
is then easy to see that there is an inclusion of
$L_1\hookrightarrow E$, which is preserved by $\ttau$, and we are
in fact the in sub-case just discussed.

Now consider the case where $\tau^* (L_1) \cong L_2 \ncong L_1$.
Then we get an isomorphism $\ttau(L_1) \cong L_2$ and
$\ttau^2(L_\nu) = L_\nu$, $\nu=1,2$. So $n$ has to be even, $E
\cong L_1\oplus L_2$ and $\ttau$ acts off diagonally. If we let
$\Xt = X/\langle\tau^2\rangle$ and $\pit : \Xt \ra Y$ be the
projection, then by applying the above, we get an equivariant
isomorphism $(\pit^*\bar E,\pit^*\bar F,w) \cong (\Lt,w_1) \oplus
(\tau^*(\Lt),w_2)$. Note in particular that $\pit^*\bar F_{\xt}$
must be either $\Lt_{\xt}$ or $\tau^*\Lt_{\xt}$ inside $\Lt_{\xt}
\oplus\tau^*\Lt_{\xt}$ for $\xt\in\pit^*(P)$. Also note that $\pit
: \Xt \ra Y$ cannot be ramified, for else $\tau : \Xt \ra \Xt$
would have fixed points, over which the action on $\Lt\oplus
\tau^*\Lt$ has eigenvalues $1$ and $-1$, contradicting the fact
that $\Lt\oplus \tau^*\Lt\conge \pit^*\bar E$. However, $\pit :
\Xt \ra Y$ is unramified if and only if $r$ is even.

Since
 $$
  \bar E = \Lt \oplus \tau^*\Lt/ \langle \tau \rangle,
 $$
we get by Lemma 2.1 in \cite{\refNR} that if $\langle \bar
L\rangle = \ker\{ \pit^* : \Pic_0(Y) \ra \Pic_0(\Xt)\}$, then
$\bar E \cong \bar E \otimes \bar L$. In fact it is easy to write
this isomorphism on $\Lt\oplus \tau^*\Lt$ explicitly as $\Phit =
\id \oplus (-\id)$ which commutes with the $\tau$ action. But then
we see that $\Phit$ takes $\pit^*\bar F$ to it self, hence we get
that $(\bar E,\bar F,w) \cong (\bar E,\bar F,w)\otimes \bar L$.

Conversely, suppose that $r$ is even and $(\bar E,\bar F,w)\in
P_a^{s,i}$, i.e. $(\bar E,\bar F,w) \cong (\bar E,\bar F,w)\otimes
\bar L$. Then $\pit : \Xt \ra Y$ is a $2$-fold unramified cover
and $\pit^* L \conge [(\check{f}_{-1})]$, where $\check{f}_-\in
\Cal M(\Xt)$ is anti invariant under $\tau$. The isomorphism
$(\bar E,\bar F,w) \cong (\bar E,\bar F,w)\otimes \bar L$ then
lifts to a $\Phit \in \aut(\pit^*\bar E,\pit^* \bar F)$ such that
$$
   \CD
     \pit^*\bar E @>\Phit>> \pit^*\bar E \\
     @V \tau^* VV @V V -\tau^*V  \\
     \pit^*\bar E @>\Phit>> \pit^*\bar E.
   \endCD
 $$
Now $\Phit$ can be normalised such that $\Phit^2 = 1$, since
$\Phit^2 \in \pit^*\aut (\bar E, \bar F) \cong \C^*$. By the
relation $-\tau^* = \Phit \tau^* \Phit$, we conclude that $\Phit
\neq \pm \id$. Hence, $\Phit$ gives an eigen decomposition $\pit^*
\bar E \cong \Lt \oplus \tau^* \Lt$ and $\pit^* \bar F$ has to be
compatible with this direct sum decomposition, because of
$\Phit$-invariance. But then it is clear that $\tF(\bar E,\bar F,
w)\in |\M^{ss}|$. \endproof

\demo{Proof of \procref{mainsetT}} We recall that $\F : P_a \ra
|\M|$ is surjective, since $\tF : P_a' \ra \Li'$ is a bijection
and $\Pi_e : \Li' \ra |\M|$ is a surjection. Let us then prove the
statement about the fibers of $\F$. In case $n$ is odd, $\Pi_e :
\Li^s \ra |\M^s|$ is a bijection, hence so is $\F : P_a^s \ra
|\M^s|$. In case $n$ is even, $\Pi_e$ is invariant under the
$\zeta_2$ action on $\Li'$ and $\Pi_e : \Li^s/\zeta_2 \ra |\M^s|$
is a bijection. Thus, by \procref{char of deg}, we get that $\F :
\bar P_a^{s,g} \ra |\M^s|$ is a bijection.

Let us now consider the fibers of $\F$ over $|\M^{ss}|$. Assume we
have $(\bar E,\bar F, w)\in P_a^{s,i}$ and $(\bar E',\bar F',
w')\in P_a^{ss}\sqcup P_a^{s,i}$. As we saw in the proof of
\procref{char of deg}, we must have that $\tF(\bar E,\bar F, w)
\cong (L\oplus \tau^*L,\ttau)$, where $L \ncong \tau^* L$ and
$\ttau(L) = \tau^*(L)$. If $(\bar E',\bar F', w')\in P_a^{s,i}$,
then there is a line bundle $L'$ which satisfies the same as $L$.
But, then since $L', L$ and $\tau^*L$ all have the same degree, we
see that either $L'\cong L$ or $L' \cong \tau^*L$. In both cases
we get an equivariant isomorphism $(L\oplus \tau^*L,\ttau) \cong
(L'\oplus \tau^*L',\ttau')$, which gives that the numeric
invariants are the same and that there is an isomorphism $(\bar
E,\bar F,w) \cong (\bar E',\bar F',w')$. If $(\bar E',\bar F',w')
\in P_a^{ss}$, we can assume that $(\bar E',\bar F',w') \cong
\Gr(\bar E',\bar F',w')$ and so $\tF(\bar E',\bar F',w') \cong
(L_1,\ttau_1)\oplus (L_2,\ttau_2)$ is a direct sum of equivariant
line bundles. From this we get that $L\cong L_1 $ or $L\cong L_2$,
because these line bundles have the same degree. This gives a
contradiction, since $L\ncong\tau^*L$. --- Summing up, we have
that $\F : P_a^{s,i}\ra \F(P_a^{s,i}) \subset |\M^{ss}|$ is a
bijection and $\F(P_a^{ss}) = |\M^{ss}| - \F(P_a^{s,i})$.

Let us now examine the fibers of $\F : P_a^{ss} \ra |\M^{ss}|$, so
let
 $$
  \tF\big((\bar L_1, w_1)\oplus (\bar L_1, w_1)\big) \cong (L_1,\ttau_1)\oplus (L_2,\ttau_2)
 $$
and
 $$
  \tF\big((\bar L'_1, w'_1)\oplus (\bar L'_1, w'_1)\big) \cong (L'_1,\ttau'_1)\oplus (L'_2,\ttau'_2)
 $$
and suppose that $L_1 \oplus L_2 \cong L_1'\oplus L_2'$. By if
necessary re-indexing, we can assume that $L_1 \cong L_1'$ and
$L_2 \cong L_2'$. It now follows immediately that in case $n$ is
odd, $[(\bar L_1,w_1)\oplus (L_2,w_2)] \sim_o [(\bar
L'_1,w'_1)\oplus (L'_2,w'_2)]$ and in case $n$ is even, that
$[(\bar L_1,w_1)\oplus (L_2,w_2)] \sim_e [(\bar L'_1,w'_1)\oplus
(L'_2,w'_2)]$ by the very construction of $\sim_o$ and $\sim_e$.
\endproof

\subhead Example: The unramified case\endsubhead Assume that $\pi
: X \ra Y$ is unramified, i.e. $P = \emptyset$. Then $\LD =
\{(0,0)\}$ and for each $\Delta \in \Dc$, $P(\bar \Delta,w) =
\bM(\bar \Delta)$ the moduli space of semi-stable holomorphic
bundles of rank $2$ on $Y$ with determinant isomorphic to $\bar
\Delta$, so
 $$
  P_a = \cases
          \bM(\bar \Delta), &\quad\text{$n$ odd} \\
          \bM(\bar \Delta')\sqcup \bM(\bar \Delta''),&\quad\text{$n$ even,}
     \endcases
 $$
where $\{\pi^*\bar\Delta',\pi^*\bar\Delta''\}=\Dc$ when $n$ is
even. The map $\F : P_a \ra |\M|$ is just $\pi^* : P_a \ra |\M|$
and if $n$ is even, the $\zeta_2$-action is just given by
tensoring with the line bundle $\bar L = \bar L_\pi^{n/2}$, where
$\langle \bar L_\pi\rangle = \ker \{\pi^* : \Pic_0(Y) \ra
\Pic_0(X)\}$. Thus, we get that
 $$
  |\M^s| \cong \cases
          \bM^s(\bar \Delta), &\quad n\text{ odd} \\
          \bM^{s,g}(\bar \Delta')/\langle \bar L\rangle \sqcup
          \bM^{s,g}(\bar \Delta'')/\langle \bar L\rangle, &\quad n\text{ even.}
     \endcases
 $$

Let us consider the case $n$ is {\it odd\/} first. If $\deg
\Delta$ is odd, then so is $\deg \bar \Delta$ and we simply get
that $|\M| \cong \bM(\bar \Delta)$. If $\deg \Delta$ is even, then
so is $\deg \bar \Delta$. Let $\bDh$ be a square root of $\bar
\Delta$. Then the map $\Xi : \Pic_0(Y) \ra \bM(\bar \Delta)^{ss}$
given by $\Xi(\bar H) = \bH \bDh \oplus \bH^{-1} \bDh$ is
surjective and invariant under the $\zeta_2$-action on $\Pic_0(Y)$
given by $\bH\mapsto \bH^{-1}$ and $\Xi : \Pic_0(Y)/\zeta_2 \ra
\bM(\bar \Delta)^{ss}$ is a bijection. The action of $\Pic_0(Y)$
on it self restricts to an action of $\langle \bar L_\pi\rangle$
on $\Pic_0(Y)$ and if $G = \langle \bar L_\pi\rangle \ltimes
\zeta_2$, then
 $$
  \pi^* \circ \Xi : \Pic_0(Y)/G \lra |\M^{ss}|
 $$
is a bijection.

Let us now consider the case where $n$ is even. Then there are two
subcases, the first being $\deg \Delta / n$ odd. We get then
immediately that
 $$
   |\M^{ss}| \cong \bM(\bar\Delta')^{s,i}\sqcup
       \bM(\bar\Delta'')^{s,i},
 $$
i.e. by Corollary 3.6 in \cite{\refNR} a disjoint union of two
copies of the Prym variety $\check{P}$ of the double cover $\pit :
\Xt \ra X$, where $\Xt = X / \langle\tau^{n/2}\rangle$.

The other subcase is $\deg \Delta / n$ even. As before we have
bijections
 $$
   \Xi' : \Pic_0(Y)/G \lra \bM(\bar \Delta')^{ss}/\sim_o
 $$
and
 $$
   \Xi'' : \Pic_0(Y)/G \lra \bM(\bar \Delta'')^{ss}/\sim_o.
 $$
But $\sim$ gives a bijection between $\bM(\bar
\Delta')^{ss}/\sim_o$ and $\bM(\bar \Delta'')^{ss}/\sim_o$, so we
get that
 $$
  |\M^{ss}| \cong \Pic_0(Y)/G \cup
    \big(\bM(\bar\Delta')^{s,i}\sqcup
    \bM(\bar\Delta'')^{s,i}\big).
 $$

Hence, the irreducible components of $|\M^{ss}|$ are the four
Kumar-varieties $\bM(\bar\Delta')^{s,i}\sqcup
\bM(\bar\Delta'')^{s,i}$, all isomorphic to $\check{P}/\{\pm 1\}$
and the quotient $\Pic_0(Y)/G$, intersecting each of the Kumar's
in finitely many points.


\head Geometric Invariant Theory Analysis of the Morphism\endhead

\seclab{git}

In this section we will prove that $\em \: \Pa \to |M|$ is a {\it
morphism\/} of varieties using geometric invariant theory (GIT),
and we will study the structure of this morphism and of the fixed
point variety. We refer the reader to \cite{\refMeS},
\cite{\refSe}, \cite{\refMuone} and \cite{\refGr} for the details
about GIT. Here we review just the minimum to fix the notation.

Let $Z$ be a smooth complex algebraic curve of genus $g= g(Z)$.
Suppose that $\Bbb E$ is the trivial bundle of rank $p = d -
2(g-1)$ over $Z$. Fix a Hilbert polynomial $\rho_0(T) = p + 2 T$,
and assume that $d > 2(2g-1)$. Let $Q_Z =
\operatorname{Quot}_{\Bbb E/Z/\Bbb C}^{\rho_0}$ be Grothendieck's
Quot scheme, \cite{\refGr}, \cite{\refSe}, and $\Cal U_{Z}$ the
universal quotient sheaf over $Q_Z \times_{\Bbb C} Z$. Then there
is an open sub-scheme $R_Z$ of $Q_Z$, characterised by the
property that $R_Z$ is exactly the points $q \in Q_Z$ for which
$\Cal U_{Z,q}$ is locally free over $Z_q = \{q\}\times_{\Bbb C} Z$
and the homomorphism $H^0(Z_q,\Bbb E) \to H^0(Z_q,\Cal U_{Z,q})$
given by the quotient morphism is an isomorphism. Then $\Cal
U_{R_Z} \deq \Cal U_Z|_{R_Z \times_{\Bbb C}Z}$ is locally free and
the sub-scheme $R_Z$ satisfies {\it local universality\/},
\cite{\refSe, Proposition~1.III.21}.

The moduli space $M$ consisting of strong equivalence classes of
semi-stable holomorphic bundles is a {\it good quotient\/} (in the
sence of GIT) of the semi-stable part $R^{ss}_X$ with respect to
action of the group $\PGL(p)$, \cite{\refSe, p.~34}.

For the parabolic case let $P \subset Z$ be a finite set of points
and denote by $$\Pj(\Cal U_{R_Z})_P = \prod_{z\in P} \Pj(\Cal
U_{Z}|_{R_Z \times_{\Bbb C} \{z\}}).$$ Let $\bar R_Z$ be the
closed subvariety of $\Pj(\Cal U_{R_Z})_P$ whose points $\{F_z\}$
satisfy that for any pair of points $(z,z')$ in $P$, the two
projections to $R_Z$ from $\Pj(|_{R_Z \times_{\Bbb C} \{z\}})$
respectively $\Pj(\Cal U_{Z}|_{R_Z \times_{\Bbb C} \{z'\}})$
agree. This is precisely what ensures that we get a (canonical)
projection
 $$
   \Pi \: \bar R_Z \longrightarrow R_Z.
 $$
Notice that a $q \in \bar R_Z$ defines a quasi parabolic structure
on the vector bundle $\Cal U_{R_Z,\Pi(q)}$ over $Z$. We shall
write $\bar\Cal U_{\bar R_Z,q}$ for this quasi parabolic bundle
and $\bar\Cal U_{\bar R_Z}$ for the corresponding bundle over
$\bar R_Z \times_{\Bbb C}Z$. Given a set of weights $w$ over $P$,
we denote the part of $\bar R_{Z}$ which is semi-stable with
respect to $w$ by $\bar R^{ss}_{Z}(w)$, and the part whose points
has a fixed determinant $L$ by $\bar R^{ss}_{Z,L}(w)$.

The moduli space $P(Z;L,w)$ of semi-stable parabolic bundles over
$Z$ with weights $w$ and determinant $L$ is a good quotient of
$\bar R^{ss}_{Z,L}(w)$ by $\PGL(p)$, \cite{\refSe, p.~84}. When
$g(Z) \ge 2$, $\bar R^{ss}_{Z,L}(w)$ is an open subset of an
irreducible variety, \cite{\refSe, p.~84}. Thus, $P(Z;L,w)$ is
irreducible. For low genus fixing, the weights and the determinant
is not necessarily enough to ensure irreducibility. See e.g.\ the
following section.

\proclaim{Lemma \n\proclab{sheaf rep map}{Lemma}} Suppose a
set-theoretic map
 $$
   f \: \bar R_{Z,L}^{ss}(w) \longrightarrow P(X;L',w'),
 $$
is presented by a locally free sheaf $\Cal V$ over $\bar
R_{Z,L}^{ss}(w) \times_{\Bbb C} X$; i.e.\ that $f(q) = [\Cal V_q]$
for all $q \in \bar R_{Z,L}^{ss}(w)$, where $[\cdot]$ means strong
equivalence class. Then $f$ is a morphism provided $\Cal V$
satisfies local universality. Moreover, if $f$ is
$\PGL(p)$-invariant it induces a morphism
 $$
   \hat f \: P(Z;L,w) \longrightarrow P(X;L',w').
 $$
\endproclaim

Of course the above statement holds for semi-stable bundles
without parabolic structures as well.

\proof By assumption, $\Cal V$ conforms to the requirements of
local universality so for every element $q \in \bar
R_{Z,L}^{ss}(w)$ there is a neighbourhood $V_q$ of $q$ in $\bar
R_{Z,L}^{ss}(w)$ and a morphism $f_q \: V_q \to \bar R_X$ so that
$\Cal V|_{V_q \times_{\Bbb C} X} \cong \check f^*\bar\Cal U_{\bar
R_{XL'}}$, where $\check f \: V_q \times_{\Bbb C} X \to \bar R_X
\times_{\Bbb C} X$ is the induced morphism. Clearly $\im f_q
\subseteq \bar R_{X,L'}^{ss}(w')$.

Let $g_X \: \bar R_{X,L'}^{ss}(w') \to P(X;L',w')$ be the good
quotient. Then $f|_{V_q} = g_X \circ f_q \: V_q \to P(X;L',w')$ is
a composition of morphisms. That means that $f$ is a morphism
locally around every point $q \in \bar R_{Z,L}^{ss}(w)$, hence, it
is a morphism.

Now if the morphism $f \: \bar R_{Z,L}^{ss}(w) \to P(X;L',w')$ is
$\GL(p)$-invariant, it gives by universality of good quotients a
morphism $\hat f$ so that
 $$
   \xymatrix{
  \bar R_{Z,L}^{ss}(w) \ar[d]_{g_X} \ar[dr]^{f}        \\
  P(Z;L,w) \ar@{.>}[r]_{\hat f}  & P(X;L',w')              }
 $$
commutes. \endproof

By choosing $Z = X$ and $\Cal V = \bar\Cal U_X|_{\bar
R_{X,L}^{ss}} \otimes p_X^*L'$ we verify that tensoring with a
line bundle $L'$ induces a morphism $P(X;L,w) \to P(X;L\otimes
{L'}^2,w)$. Similarly, we see that pullback with respect to
morphisms between curves induces morphism between corresponding
moduli spaces. In particular we get that the action of
$\langle\tau \rangle$ on $M$ is algebraic.\pagelab{pullback morph}

\proclaim{Proposition \n\proclab{mod sp morphism}{Proposition}}
The map
 $
   \em \: P_a \to |M|
 $
is a projective morphism of varieties.\hfill$\square$\endproclaim

Note that we in this propostion really are considering the
restriction of this map to each irreducible component of $P_a$
mapping to the corresponding target component in $|M|$.

\proof Let $\Delta\in\Dc$ and consider the corresponding
$\bar\Delta \in \Pic_0(Y)$. By altering $\bar \Delta$ (and
therefore also $\Delta$ correspondingly) by a sufficient high
power of an ample line bundle over $Y$, we may assume we are in
the realm of the GIT-construction of $P(\bar \Delta,w)$ and
$M(\Delta)$. Formula (2.24) gives an explicit expression for the
lift $\bar \em : \bar R^{ss}_{Y,\bar\Delta} \ra |M(\Delta)|$ of
$\em : P(\bar \Delta,w) \ra |M(\Delta)|$. Tensoring by $[D_2]$
clearly induces morphisms as remarked before, so we only need to
concentrate on the iterated elementary modification. In order to
use \procref{sheaf rep map}, we need to construct a locally free
sheaf $\Cal V$ representing these iterated elementary
modifications $\Gamma_{(wn',\cdot)}$. Let $\Cal G_0$ be the sheaf
over $\bar R_{Y,\bar \Delta}^{ss}(w) \times_{\Bbb C} X$ defined by
the flag such that $\Cal G_{0,((\bar E, \bar F, w), x)} =\pi^*\bar
F_x \subseteq \pi^*\bar E_x$ for $x \in \pi^{-1}(P_w)$, $P_w =
\{\, y \in P \mid w(y) \ne 0 \, \}$, and $0$ elsewhere. Put $\Cal
V_0 = (\id_{\bar R_{Y,\bar \Delta}(w)} \times \pi)^*\bar\Cal
U_{\bar R_{Y,\bar \Delta}(w)}$, and consider the ``universal
skyscraper sheaf'' $\Cal S_0$ with stalks
 $$
   \Cal S_{0, ((\bar E, \bar F, a), x)} = \cases
        \Cal V_{0, ((\bar E, \bar F, w), x)} /
                \Cal G_{0, ((\bar E, \bar F, w), x)},
                        \quad \text{for } x \in \pi^{-1} (P_w) \\
        0, \quad \quad\quad\quad\quad \quad \quad\quad\quad\quad\text{elsewhere.} \endcases
 $$

There is a natural surjection
 $$
   \lambda_0 \: \Cal V_0 \longrightarrow \Cal S_0 \to 0.
 $$
Put $\Cal V_1 = \ker\lambda_0$, let and $\iota_0 \: \Cal V_1 \to
\Cal V_0$ be the canonical sheaf inclusion. As in \secref{prelim},
$\ker \lambda_0$ is locally free, the natural lift $\tilde\tau$ to
$\Cal V_0$ induces a lift $\tilde\tau$ to $\Cal V_1$ and this
gives a sheaf $\Cal G_1$ defined as the eigenspaces for
$\tilde\tau_x$, $x \in \pi^{-1} (P_w)$, that are not annihilated
by the morphism on fibers, induced by $\iota_0$. Then define $\Cal
S_1$, $\lambda_1 \: \Cal V_1 \twoheadrightarrow \Cal S_1$, and
$\Cal V_2 = \ker\lambda_1$ in the same way and continue the
process through to $\Cal S_M$, $M$ being the number of steps
necessary to do $\elm_{(wn',\cdot)}$, and define
 $$
   \Cal V = \Cal V_{M+1} = \ker\lambda_M.
 $$
By the alteration of $\bar\Delta$, the degree is constrained such
that $\Cal V$ satisfies local universality.

The morphism represented by $\Cal V$ is $\GL(\bar p)$-invariant
(because $\bar\Cal U_{\bar R_{Y,\bar \Delta}(w),q} \cong \bar\Cal
U_{\bar R_{Y,\bar \Delta}(w),A\cdot q}$ for all $A \in \GL(\bar
p)$). Therefore, there is an induced morphism from the moduli
space $P(\bar \Delta,w)$ to $|M(\Delta)|$ represented by $\Cal V$.
But clearly this morphism is $\em$.

That $\em$ is actually projective in the sense of Grothendieck,
follows immediately since $|M|$ is separated over $\Spec\Bbb C$
and both $P_a$ and $|M|$ are projective (over $\Spec\Bbb C$).
\endproof

Our main statement about the algebraic geometric properties of
$\em$ reads as follows:

\proclaim{Theorem \n\proclab{em normal even}{Theorem}\proclab{em
normal odd}{Theorem}} When $n$ is odd, the morphism
 $$
   \em \: P_a \longrightarrow |M|,
 $$
is a birational equivalence which is the normalisation morphism of
each of the irreducible components of $|M|$. When $n$ is even, the
morphism $\em$ factors through the $\zeta_2$-quotient and gives a
projective morphism
 $$
   \em \: P_a\big/\zeta_2 \longrightarrow |M|,
 $$
which is a birational equivalence. When restricted to each of the
irreducible components of $P_a/\zeta_2$ it is likewise the
normalisation morphism for the corresponding irreducible component
of $|M|$. \endproclaim

\proof It follows from \procref{act2} and \noprocref{sheaf rep
map} and arguments similar to what we saw in the proof of
\procref{mod sp morphism}, that the $\zeta_2$-action on $P_a$ is
algebraic. General theory gives that the quotient of a projective
variety under the action of a finite group acting algebraically is
again a projective variety. Moreover, if the variety is normal,
then so is the quotient. Hence, in the case $n$ is even, we see
that $P_a/\zeta_2$ is a projective variety whose irreducible
components are normal. It is a consequence of Zariski's Main
Theorem, \cite{\refGrfour, Corollary (4.4.9)}, that the
restriction of $\F \: P_a/\zeta_2 \to |M|$ to each of these
components is the normalisation morphisms. This is because it is
projective and therefore proper, and since it is generically a
bijection, it is both a birational equivalence and a dominant
morphism. The result then follows from the universal properties of
the normalisation of a variety. The same arguments applies to $\F
\: P_a \to |M|$ when $n$ is odd. \endproof

\remark{Remark \n\proclab{when normal}{Remark}} We may split up
the equivalence relations $\sim_o$ resp.\ $\sim_e$ into an
equivalence relation $\sim_n$ acting only within each component,
and an equivalence relation acting strictly between components. We
remark that a component of $|M|$ is itself normal if and only
$\sim_n$ is trivial on the corresponding component of $P_a$ resp.\
$P_a/\zeta_2$.\endremark

\head The Hyperelliptic Involution\endhead \NObreak
\seclab{hypell} \NObreak

Let $X$ be a compact, hyperelliptic curve of genus $g \ge 2$ with
a hyperelliptic involution $J$. Then $\pi \: X \to X/\langle
J\rangle = \Pj^1$ and the set of fixed points of $J$ in $X$ are
exactly the Weierstrass points $x_1, \dots x_{2g + 2}$. We denote
$\pi(x_j) = z_j$ and let $P = \{z_1, \ldots z_{2g+2}\}$. So in
this case $n=2$ and for each $y\in P$ we see that $k=1$, $n' = 2$
and $\xi = -1$. Let us now apply \procref{mainsetT} to describe
the fixed variety $|M|$ of $J$ in $M$, the moduli space of
semi-stable holomorphic bundles with trivial determinant.

In this case $\Dc = \{\Delta_0,\Delta_1\}$, where $\Delta_0 = 0 $
and $\Delta_1 = (f_-)$ with $f_-$ a meromorphic function such that
$f_- = - f_- \circ J$. If we choose an identification of $\Pj^1 =
\C \cup \infty$, such that $\infty \not\in P$, then we have
explicitly that $f_-$ is the anti-invariant meromorphic function
on $X$ determined by the multivalued meromorphic function $$ z
\longmapsto \sqrt{\prod_{i=1}^{2g+2}\big(z-z_i\big)} $$ on
$\Pj^1$, whose associated ramified cover exactly is $\pi : X \ra
\Pj^1$. Thus,
 $$
  \Delta_1 = \sum_{i=1}^{2g+2}x_i -(g+1)\cdot \pi^{-1}(\infty).
 $$

Let us examine the case $\Delta_0\in \Dc$ first. In this case $$
  \Lambda_{\Delta_0} = \{(0,0),(1,1)\}^{\times (2g+2)}.
$$ Hence, $d_1(y) = d_2(y)$ so $w(y) = 0$ for all $y\in P$. An
element of $\Lambda_{\Delta_0}$ is just determined by the subset
$Q \subseteq P$ given by $Q = \{\,y\in P \mid d_1(y) = d_2(y) =1\,
\}$ and the corresponding line bundle on $\Pj^1$ given by formula
\tagref{bardelta} is $$
  \bar \Delta_Q = [-Q] = \so(-d_Q)
$$ where $d_Q=|Q|$. Since we only allow quasi-parabolic structures
at the points, where the weights $w$ are non-zero, there are no
quasi-parabolic structures to consider in this case and $$
P_a|_{\Delta_0} = \coprod_{Q \subseteq P} M\big(\so(-d_Q)\big), $$
where $M\big(\so(-d_Q)\big)$ is the moduli space of semi-stable
bundles on $\Pj^1$ with determinant $\so(-d_Q)$. Grothendieck's
classification of vector bundles on $\Pj^1$ combined with
semi-stability implies that $$
  M\big(\so(-d_Q)\big) = \cases \emptyset &\quad\text{if $d_Q$ is
            odd},  \\
        \big\{\so\big(-\tfrac{d_Q}2\big) \oplus
                \so\big(-\tfrac{d_Q}2\big)\big\} &\quad\text{if $d_Q$ is even}
        .\endcases
$$ Denote the point in $P_a$ corresponding to $Q \subseteq P$ by
$\llangle Q\rrangle$, then
 $$
  P_a|_{\Delta_0} = \coprod\Sb Q \subseteq P \\ d_Q \text{ even}\endSb
        \llangle Q\rrangle.
 $$

Let us now examine the case $\Delta_1\in \Dc$. Then
$\Lambda_{\Delta_1}$ consist only of one element, namely $$
  \Lambda_{\Delta_1} = \{(0,1)\}^{\times (2g+2)},
$$ and we have that $w(y) = \frac12$ for all $y\in P$. We see that
$\bar \Delta_1 = \so\big(-(g+1)\big)$, so in this case $$
P_a|_{\Delta_1} = P\left(\so\big(-(g+1)\big),\tfrac12\right), $$
where $\frac12$ means weight $\frac12$ at all points of $P$.

Let $d = -(g+1)$ be the degree of the underlying bundle and let us
further analyse this parabolic moduli space $\Phf =
P\left(\so(d),\tfrac12\right) = P_a|_{\Delta_1}$. If $(\bar E,\bar
F,\frac12) \in \Phf$, then there is an integer $c$ such that $\bar
E \cong \so(c) \oplus \so(d-c)$. We see that $\pa\mu(\bar E) = 0$.
By symmetry, we can assume that $d-c \le \frac d2 \le c$.

If $\bar L$ is a parabolic subbundle of $\bar E$ then
$\pa\deg(\bar L) = \deg(\bar L) +g +1 -\frac12 |P(\bar L)|$, where
$P(\bar L) = \{\, y \in P \mid \bar L_y \ne F_y \, \} $, so the
semi-stability condition reads
 $$
   \deg(\bar L) +g+1 \le \frac{|P(\bar L)|}2
 $$
with ``$<$'' substituted for ``$\le$'' if we want stability. In
particular, we notice that $\deg(\bar L) \leq 0$, since $ |P(\bar
L)| \le 2(g+1)$. If we apply this to the subbundle $\so(c)$, we
get that $c \le 0$. Let now $\Phf^s_c$ be the moduli space of
stable parabolic bundles over $\Pj^1$ with parabolic weights
$\frac12$ over each point in $P$, such that the underlying vector
bundle is isomorphic to $\so(c) \oplus \so(d-c)$. Then clearly
$\Phf^s_0 = \emptyset$ and we have that
 $$
   \Phf^s = \coprod_{\frac{d}2 \leq c < 0}
   \Phf^s_c.
 $$
The subsets $\Phf_c^s$ are of course disjoint, however for $c > 0$
their closure in $\Phf$, which we will denote $\Phf_c$, do
intersect in the subsets represented by only semi-stable bundles,
as we shall now see. For each subset $Q\subseteq P$, such that
$d_Q$ is even, there is a semi-stable parabolic structure on
$\so\big(\frac{-d_Q}2\big) \oplus \so\big(d+\frac{d_Q}2\big)$ with
all weights $\frac12$, namely over points in $Q$, let the flag be
the fiber of the first summand and over the rest of the points,
choose the fiber of the second summand. Let us denote this
parabolic bundle $\langle Q \rangle$. It is clearly semi-stable,
since it is the direct sum of two semi-stable parabolic bundles
with parabolic slope $0$. We define $\Phf_0$ to be the subset of
$\Phf$ whose underlying bundle can be represented by $\so(0)
\oplus \so(d)$. Thus, $\Phf_0 = \{\,\langle \emptyset\rangle\,\}$.

\proclaim{Proposition \n\proclab{parabolhalfP}{Proposition}} When
$ c > 0$, each of the moduli spaces $\Phf^s_c$ is a non-empty
connected quasi-projective variety whose closure $\Phf_c$
 is an irreducible component of $\Phf$.
The set of S-equivalence classes of semi-stable, but not stable,
parabolic bundles in $\Phf$ are in bijective correspondence with
the even cardinality subsets of $P$ by the above construction.
Moreover, the closure $\Phf_c$ is explicitly
 $$
   \Phf_c = \Phf^s_c \sqcup
   \coprod\Sb Q\subseteq P \\
        \text{even } d_Q \leq -2c\endSb
      \langle Q\rangle.
 $$
\endproclaim

We see in particular from this proposition that the intersection
of all the irreducible components consists of exactly one element
 $$
   \langle\emptyset \rangle = \bigcap_{\frac{d}2 \leq c \leq 0}
    \Phf_c = \Phf_0
 $$
and hence, that $\Phf$ is in fact connected.

\proof First we shall establish that for each $\frac{d}2 \le c <
0$, there exist stable parabolic vector bundles whose underlying
vector bundle is $\so(c) \oplus \so(d-c)$. The possible set of
flags $F$ is parameterised by $(\Pj^1)^{2g+2}$ corresponding to
the $2g+2$ parabolic points. Stability means that for every
inclusion $i : \so(e) \ra \so(c)\oplus \so(d-c)$, we must have
that $e + g+1 < \frac12 \big|P\big(i(\so(e))\big)\big|$. We see
immediately that for such an inclusion to exist, we must have
$e\leq c$, and if $d-c < e \leq c$, then $e=c$ and the inclusion
is forced to be $i : \so(c) \ra \so(c)\oplus \{0\}$. For $e < d$
the stability condition is trivially satisfied. Hence, the only
interesting range for $e$ is $d \leq e \leq d-c$. For a given $e$
in this range let $I^c_e$ be the set of inclusions of $\so(e)$
into $\so(c)\oplus \so(d-c)$, modulo automorphisms of $\so(e)$. We
note that $I^c_e$ is an open dense subset of
 $$
   \Pj\Big(\Hom\big(\so(e),\so(c)\big)\oplus
            \Hom\big(\so(e),\so(d-c)\big)\Big)
        \cong\Pj^{-g-2e},
 $$
each of which are linearly mapped into the set of flag
$(\Pj^1)^{2g+2}$. Hence, we see that $\cup_{d \leq e \leq
d-c}I^c_e$ has codimension $g$. Since $g\geq 1$, this means there
are always flag configurations, which are not contained in the
subset $\cup_{d \leq e \leq d-c}I_e^c$. Therefore, the set of flag
configurations, which makes $\so(c)\oplus \so(d-c)$ stable is a
non-empty, open, dense subset of $(\Pj^1)^{2g+2}$. Since
$\Phf^s_c$ is the $\Aut\big(\so(c)\oplus \so(d-c)\big)$-quotient
of this non-empty dense subset, we get the desired conclusion
about the parabolic moduli spaces $\Phf^s_c$.

Suppose we have two parabolic line bundles $(L_i,w_i)$, $i=1,2$,
of zero parabolic degree, such that $(L_1,w_1)\oplus(L_2,w_2)$
represents the S-equivalence class of a semi-stable parabolic
bundle in $\Phf$. Then we can assume there is an $\frac{d}2 \leq
e\leq 0$ integer such that $L_1 \cong \so(e)$ and $L_2 \cong
\so(d-e)$. Since $(L_i,w_i)$, $i=1,2$, have zero parabolic degree,
there is a subset $Q\subseteq P$ of cardinality $2e$, such that
$z\in Q$ if and only if $w_1(z) = \frac12$ and $z\in P-Q$ if and
only if $w_2(z) = \frac12$. Hence, $(L_1,w_1)\oplus(L_2,w_2) \cong
\langle Q\rangle$. In order for $\langle Q\rangle$ to be contained
in $\Phf_c$ we just need that $\so\big(d+\frac{d_Q}2\big)$ can be
included in $\so(c)\oplus\so(d-c)$, which is the case if an only
$d+\frac{d_Q}2 \leq d-c$, i.e. if and only if $ d_Q \leq -2c $.
\endproof

By analyzing how the automorphism group of the bundle
$\so(c)\oplus\so(d-c)$ acts on the flags of the bundle, we arrive
at the following dimension formula $$
  \dim \Phf_c = \cases g - 2c -1,
            &\quad\text{if } 0 > c > \frac d2 = -\frac{g+1}2 \\
      2g -1, &\quad\text{if } c = \frac d2 = -\frac{g+1}2.
  \endcases \tagn
$$\numlab{en label tak}

Now, let us consider the $\zeta_2$-action and equivalence
relations $\sim_o$ and $\sim$
 on $P_a$:

A moments examining of \procref{DEFeqe} and \procref{sime} gives
that $\sim$ bijectively relates $P_a^{ss}|_ {\Delta_0} =
\coprod_{Q \subseteq P, d_Q \text{ even}} \llangle Q\rrangle$ and
$P_a^{ss}|_ {\Delta_1} = \coprod_{Q \subseteq P, d_Q \text{ even}}
\langle Q\rangle.$ Namely $\llangle Q_1\rrangle \sim \langle
Q_2\rangle$ if and only if $Q_1 = Q_2$. Since $P_a|_ {\Delta_0} =
P_a^{ss}|_ {\Delta_0}$, we see that $P_a|_ {\Delta_0}$ gets
identified with $P_a^{ss}|_ {\Delta_1}$. As $n=2$ the
$\zeta_2$-action coincides with the relation $\sim_o$ on
$P_a^{ss}$, the induced relation on $P_a^{ss}/\zeta_2$ is trivial
from which it follows that the components of $|M|$ will be normal.
In fact, it is not difficult to see that the action is trivial on
$P_a^{ss}$. Since $\sim_o$ and $\sim$ are compatible, we obtain
that
 $$
   P_a / (\sim_e, \zeta_2) = \Phf/\zeta_2,
 $$
where we note that the $\zeta_2$-action preserves each of the
irreducible components $\Phf_c$.

We conclude:

\proclaim{Proposition \n} Let $X$ be a hyperelliptic curve of
genus $g \ge 1$ and let $J$ be a hyperelliptic involution. The
fixed point set $|M|$ is decomposed into irreducible components
 $$
   |M| =  \Phf/\zeta_2 = \bigcup_{\frac{d}2 \le c < 0}
        \Phf_c\big/\zeta_2.
 $$
The irreducible components $\Phf_c/\zeta_2$, $c > 0$, are normal
sub-varieties of $|M|$ of the dimensions stated in formula
\tagref{en label tak}. They only intersect within the finite
subset $|M^{ss}|$, and their intersections are given by
\procref{parabolhalfP}. \hfill$\square$ \endproclaim

We note that the component of $|M|$ of maximal dimension is
$\Phf_c$, with $c = -\frac{g+1}2$ in case $g$ is odd and $c=
-\frac{g}2$ in case $g$ is even. The dimension of this component
is $2g-1$.

\appendix The Kernel of $\pi^* :$ {\rm Pic}$(Y) \to $ {\rm
Pic}$(X)$\endappendix A

\seclab{kernel sec}

In this section we calculate the kernel $\ker\{\, \pi^* \: \Pic(Y)
\to \Pic(X)\, \}$ for the covering map $\pi \: X \to Y = X/
\langle\tau\rangle$. We will refer to this kernel as $\ker\pi^*$.

Assume that $D \in \Div(Y)$ is a divisor for which the associated
line bundle pulls back to the trivial bundle, $\so_X$, on $X$,
i.e.\ $\pi^*D = (f)$ for some meromorphic function $f \in \Cal
M(X)$. Then using the relation between the pullback $\pi^*$ and
the norm map $\Nm$,
 $$
   n \cdot D = \Nm \circ \pi^*(D) = \Nm\big((f)\big) = \big(\Nm(f)\big),
 $$
we see that every element of the kernel has order at most $n$;
hence, the kernel is a subgroup of the $n$-torsion points
$\Pic_0^{(n)}(Y)$ in $\Pic(Y)$. The analysis of which subgroup, is
divided into two cases.

Before we proceed we shall need to prove the following general
statement:

\proclaim{Lemma \n\proclab{trans mero}{Lemma}} Let $\tau \: X \to
X$ be an automorphism of order $n$ of a curve $X$. Then for any
$n$'th root of unity, $\mu$, there exists a meromorphic function
$h \in \Cal M(X)$ such that
 $$
   h \circ \tau = \mu \cdot h.
 $$
\endproclaim

\proof The space of meromorphic functions $\Cal M(Y)$ on $Y$ is a
subfield of $\Cal M(X)$. In fact as $\pi^*\Cal M(Y) = \Cal
M(X)^{\langle \tau\rangle}$ are the $\tau$-invariant meromorphic
functions, it follows from a theorem of Artin, \cite{\refL,
Theorem~VIII.1.8}, that $\Cal M(X)$ is a cyclic Galois extension
of $\Cal M(Y)$ of degree $n$. It is then a consequence of
Hilbert's Theorem~90, \cite{\refL, Theorem~VIII.6.1}, that there
is a function $h \in \Cal M(X) - \{0\}$ such that the element
$\varepsilon \in \zeta_n \subset \Cal M(X)$ satisfies that
$\varepsilon = h/(h\circ\tau)$. I.e.\ $h \circ \tau = \mu\cdot h$,
if $\mu = \varepsilon^{-1}$. \endproof

Using \procref{cover factor} we split up the calculation of the
kernel into two; one regarding the unramified case and the other
the {\it completely\/} ramified case:

\proclaim{Lemma \n\proclab{unram-ker}{Lemma}} Let $\tau \: X \to
X$ be an automorphism of order $n$ without any special orbits, so
that the covering projection $\pi \: X \to Y =
X/\langle\tau\rangle$ is unramified. Then there is a line bundle
$\bar L_\pi \in \Pic(Y)$ of order $n$ such that
 $$
   \ker\pi^* = \langle \bar L_\pi \rangle.
 $$
\endproclaim

\proof Let $\mu \in \zeta_n$ be a prime root of unity and $h \in
\Cal M(X)$ the function from \procref{trans mero} such that $h
\circ \tau = \mu \cdot h$. Then the norm map of $h$ is given by
 $$
   \Nm(h)(y)
        = \prod_{i = 0}^{n-1} h\big(\tau^i(x)\big)
        = \mu^{\sum_{i= 0}^{n-1} i}\cdot h(x)^n
        = \mu^{\frac{n(n-1)}{2}}\cdot h(x)^n
        = (-1)^{n-1} h(x)^n,
 $$
for any $x \in \pi^{-1}(y)$. We see that the divisor
$\big(\Nm(h)\big)$ is divisible by $n$. Notice also that the
pullback is given by $\pi^*\Nm(h)(x) = \Nm(h)\big(\pi(x)\big) =
(-1)^{n-1} h(x)^n$. Let $D \in \Div(Y)$ be such that $n D =
\big(\Nm(h)\big)$, then
 $$
   n\pi^*D = \big(\pi^*\Nm(h)\big) = \big((-1)^{n-1}h^n\big) = (h^n) = n(h),
 $$
and thus $\pi^*D$ is a principal divisor with $\pi^*D = (h)$;
i.e.\ $[D] \in \ker\pi^*$.

Suppose that $E \in \Div(Y)$ with $[E] \in \ker\pi^*$. Then
$\pi^*E = (f)$ for some $f \in \Cal M(X)$. As $\pi^*E$ is
$\tau$-invariant, there is a root of unity $\mu' = \mu^j$, $0 \le
j \le n-1$, so that $f \circ \tau = \mu' \cdot f$, and $ \pi^*(jD
- E) = \left(\frac{h^j}f\right)$. But
 $$
   \frac{h^j}f \circ \tau =  \frac{\mu' \cdot h^j}
                {\mu'\cdot f} = \frac{h^j}f,
 $$
so there must be a $g \in \Cal M(Y)$ with $\pi^*g = \frac{h^j}f$,
and in that case
 $$
   jD - E = (g).
 $$
This means that $\ker\pi^* = \langle [D] \rangle$.

We know that $[D]$ has order at most $n$ since $nD =
\big(\Nm(h)\big)$, so suppose that $jD = (g)$ for some $g \in \Cal
M(Y)$. Then $(\pi^*g) = (h^j)$ in which case there must exist a
$\mu' \in \so^*(X) = \Bbb C^*$ such that $\pi^*g = \mu'\cdot h^j$.
But $\pi^*g$ is $\tau$-invariant while $h^j \circ \tau = \mu^j
\cdot h^j$, so $\pi^*g = \mu'\cdot h^j$ if and only if $j = 0 \mod
n$. This means $\bar L_\pi = [D]$ is indeed of order $n$.
\endproof

Notice that the meromorphic function $h$ in the proof is not
unique, since we can modify it by the pullback of any meromorphic
function on $Y$. Hence, $D$ isn't unique either, but there is a
unique class of linearly equivalent divisors, and thus, $\bar
L_\pi$ is in fact unique given $\mu\in\zeta_n$.

As a final remark to this we observe that if $n = n_1 \cdot n_2$
then an unramified $\pi$ factors through $\pi_1 \: X \to X_1 = X/
\langle\tau^{n_2}\rangle$ and the induced $\pi_2 \: X_1 \to Y$.
Then there are divisors $D$, $D_2 \in \Div(Y)$ and $D_1 \in
\Div(X_1)$ as above, such that $\ker\pi^* = \langle[D]\rangle$,
$\ker\pi_2^* = \langle[D_2]\rangle$ and $\ker\pi_1^* =
\langle[D_1]\rangle$, where $D_2 = n_1 D$ and $\pi_2^*D = D_1$.
This is because the meromorphic functions $h$, $h_1 \in \Cal M(X)$
and $h_2 \in \Cal M(X_1)$ can be chosen as $h_1 = h$ and $h_2 =
\Nm_{\pi_1}(h)$.

\proclaim{Lemma \n\proclab{prime-ker}{Lemma}} Let $\tau \: X \to
X$ be an automorphism such that the projection $\pi \: X \to Y =
X/\langle\tau\rangle$ is a ramified covering and such that the
orbit lengths $\{\, k(y) \mid y \in Y \, \}$ are co-prime. Then
the pullback
 $$
   \pi^* \: \Pic(Y) \longrightarrow \Pic(X)
 $$
is injective.\endproclaim

\proof The proof of this statement is essentially due to Mumford
in \cite{\refMutwo}. Let $q \: L \to Y$ be a line bundle which is
a $n$-torsion point in $\Pic(Y)$, and define the curve
 $$
   X_L \deq \{\, s \in L \mid s^n = 1 \in L^n \cong \so_Y \, \}
 $$
and the unramified covering
 $$
   \pi_L \: X_L \longrightarrow  Y.
 $$

Notice that a global non-zero section $s$ of $L$, if it exits, can
be scaled to satisfy $s^n = 1$, so that it gives a global section
of $X_L$. On the other hand, a global section of $X_L$ is a global
section of $L$. Hence, $L$ is a trivial bundle if and only if
$X_L$ is a trivial covering. Since
 $$
  \align
   \pi^* X_L &= \{\, (x,\hat x) \in X \times X_L \mid
                        \pi(x) = \pi_L(\hat x) \, \}
        = \{\, \hat s \in \pi^*L \mid \hat s^n = 1 \, \}
        = X_{\pi^*L}
  \endalign
 $$
it follows by the same argument that $\pi^*L$ is a trivial line
bundle if and only if $\pi^*X_L$ is a trivial covering of $X$.

Suppose now that $L \in \ker\pi^*$ so that $\pi^*X_L$ is trivial
and $\varphi \: X \times \{\, 1, \dots n\, \} \to \pi^*X_L$ is a
trivialisation. That gives rise to a commutative diagram
 $$
   \xymatrix{
     X \times [n] \ar[d]_{\pr} \ar[r]^{\varphi}_{\cong}
            & \pi^*X_L \ar[d]_{\pr_1} \ar[r]^{\pr_2} & X_l
         \ar[d]^{\pi_L} \\
     X \ar@{=}[r] & X \ar[r]^{\pi} & Y   }
 $$
where $[n] = \{\, 1, \dots, n\, \}$. The projection $\pr$ has
obvious sections $\sigma \: X \to X \times [n]$, and the
composition $\psi = \pr_2 \circ \varphi \circ \sigma$ is a
morphism of coverings
 $$
   \xymatrix{
  X \ar[rr]^{\psi} \ar[dr]_{\pi}
                &  &    X_L \ar[dl]^{\pi_L}    \\
                & Y                 }
 $$

Now for every $x \in X$ there is an $n$'th root of unity $\mu(x)$
so that $\psi\circ\tau(x) = \mu(x)\cdot\psi(x)$. If $k(x) = k\circ
\pi(x)$ denotes the length of the orbit through $x$,
$\tau^{k(x)}(x) = x$ so $\mu^{k(x)} = 1$ for all $x$ which means
that $\ord\mu(x) \mid k(x)$. By the basic assumption that the
action of $\tau$ does not split into components, it follows that
for each of the finitely many $k(x) \in \{\, k(y) \mid y \in Y \,
\}$ and any component $X_\alpha$ of $X$ there is an $x' \in
X_\alpha$ through which there is an orbit of length $k(x')= k(x)$.
By continuity, $\mu$ is constant on each component of $X$ and thus
$\ord \mu \mid k(x)$ for all $x$, so $\ord\mu \mid \gcd\{\, k(y)
\mid y \in Y \, \} = 1$. Hence, $\mu = 1$ which in turn means that
$|\im\psi \cap \pi_L^{-1}(x)| = 1$ for every $x \in X$. But
 $$
   n = \deg\pi = \deg(\pi_L|_{\im\psi}) \cdot \deg\psi =
   \deg(\pi_L|_{\im \psi})\cdot n,
 $$
so $\deg(\pi_L|_{\im\psi}) = 1$ and
 $$
   (\pi_L|_{\im\psi})^{-1} \: Y \longrightarrow \im\psi \subseteq
   X_L
 $$
is a global section. Hence, $L \cong \so_Y$.\endproof

Combining \procref{unram-ker} and \procref{prime-ker} we get the
general statement:

\proclaim{Proposition \n\proclab{ker}{Proposition}} Let $X$ be a
smooth algebraic curve and $\tau \: X \to X$ an automorphism of
order $n$ with possible special orbits and let $\pi \: X \to Y =
X/ \langle\tau\rangle$ be the induced covering. Then there is a
line bundle $\bar L_\pi$\pagelab{kerpi gen}\ over $Y$ such that
 $$
   \ker\pi^* = \langle \bar L_\pi\rangle,
 $$
and the order of $\bar L_\pi$ is the greatest common divisor, $r =
\gcd\{\, k(y) \mid y \in Y \, \}$, of the orbit lengths. \hfill
$\square$ \endproclaim


\references References\endreferences

 \enddocument